\documentclass[%
 reprint,
 amsmath,amssymb,
 aps,
]{revtex4-1}

\usepackage{graphicx}
\usepackage{dcolumn}
\usepackage{bm}
\usepackage{natbib}
\begin{document}

\preprint{APS/123-QED}

\title{Synchronization of Oscillators via Active Media}

\author{Derek Orr}
 \email{djo15@pitt.edu}
\author{Bard Ermentrout}%
 \email{bard@pitt.edu}
\affiliation{%
 University of Pittsburgh\\
}%




\date{\today}

\begin{abstract}
In this paper, we study pairs of oscillators that are indirectly coupled via active (excitable) cells. We introduce a scalar phase model for coupled oscillators and excitable cells. We first show that one excitable and one oscillatory cell will exhibit phase locking at a variety of $m:n$ patterns. We next introduce a second oscillatory cell and show that the only attractor is synchrony between the oscillators. We will also study the robustness to heterogeneity when the excitable cell fires or is quiescent. We next examine the dynamics when the oscillators are coupled via two excitable cells. In this case, the dynamics are very complicated with many forms of bistability and, in some cases, chaotic behavior.  We also apply weak coupling analysis to this case and explain some of the degeneracies observed in the bifurcation diagram. Further, we look at pairs of oscillators coupled via long chains of excitable cells and show that small differences in the frequency of the oscillators makes their locking more robust. Finally, we demonstrate many of the same phenomena seen in the phase model for a gap-junction coupled system of Morris-Lecar neurons.
\end{abstract}

\maketitle


\section{\label{sec:level1}Introduction}

How coupled oscillators synchronize is an important and much studied phenomenon. In many analyses of coupled oscillators, each element is a limit cycle oscillator and techniques such as weak coupling are applied. However, many systems, such as neurons, have conditional oscillators; that is, they oscillate only when given enough drive. This is the idea of an excitable system that has a unique globally attracting equilibrium point, but with a sufficiently large perturbation, it can oscillate once before returning to the stable equilibrium. What happens when there are oscillators coupled indirectly via excitable systems remains an open question. For example, in the early stages of aggregation of cellular slime molds, each cell is excitable, but some of the cells become oscillatory and the result is a global oscillatory system that induces the organism to ultimately organize into a slug (see \cite{dicty}). Within the smooth muscle of the intestine are a small number of spontaneously active cells (interstitial cells of Cajal) that are coupled and organized to form waves through the intervening non-oscillatory cells (see \cite{gut1,gut2}). Interactions between oscillatory and non-oscillatory glial cells are thought to underlie synchronization for circadian rhythms (see \cite{glia}).  

This problem has been studied in the context of global all-to-all coupling where each element is connected to all the other elements.  For example, in \cite{Daido1} they couple two populations of Stuart-Landau equations where one of the populations has a stable equilibrium and the other is oscillatory and analyze the ensemble dynamics as the relative numbers of active and inactive elements varies.  Similarly, in \cite{DePazo} the authors analyze sinusoidally coupled mixtures of oscillators and excitable cells where each cell is represented as a scalar phase model. Like \cite{Daido2}, they study the onset of collective synchrony as the ratio of oscillators changes. Others have used the so-called Ott-Antonsen reduction to study collective dynamics of mixtures of oscillatory and excitable elements (see \cite{Stro,Luke2,Luke1}).

The previous work on this problem relies on the fact that all elements are globally coupled to each other. On the other hand, the biological examples we have described are much more locally coupled. This is the scenario that we focus on in this paper. Synchrony between two neurons coupled via passive dendrites has been studied in \cite{crook}. Others have broadened this to include ``quasi-active" dendrites, though it is still a linear theory (see \cite{svenson-coombes,goldberg}). As an initial attempt to understand interactions between oscillators and excitable systems, we are interested in describing the dynamics between two oscillatory cells distributed in a simple chain with intervening excitable cells between them.  We will use a simple one-dimensional model for both the excitable and the oscillatory cells and then show that similar phenomena hold in more realistic neural models. 

We first introduce the class of models that will be our focus and then analyze small chains where there are one or two excitable cells between the two oscillators.  We will vary coupling strengths, degree of excitability, and heterogeneities in the oscillators.  In the case where the effects of the oscillators on the excitable cells are sufficiently small that they cannot induce the excitable cells to fire, we apply weak-coupling analysis and show that the results match the behavior of the full system.  We demonstrate a variety of different locking regimes as well as complex chaotic behavior. We also briefly look at longer chains of excitable cells and show that small differences in the oscillators make locking between them much more likely when the excitable chains are long.  We finally show that similar dynamics in a gap-junction coupled biophysical model and conclude with a discussion about future directions.

\section{\label{sec:level1}Methods}

There are two broad types of excitability \cite{rinz-erm89}: Class II which occurs for a system near a sub-critical Hopf bifurcation and Class I, which occurs when there is a saddle-node infinite cycle (SNIC) bifurcation. The latter type of excitability lends itself to simple one-dimensional dynamics on a circle \cite{shinimoto,erm81}, thus this will be the type of excitability we will consider in the paper.  The simplest version of this excitability takes the form
\begin{equation}
\label{eq:exc}
\frac{dy}{dt}=1-b\cos(y)=:f(y)
\end{equation} 
where $b\ge0$ is a parameter and $y\in [0,2\pi)$ lies on the circle.  When $b>1$, then (\ref{eq:exc}) has two equilibria $y^{\pm}=\pm \arccos(1/b)$, with $y^-$ (``rest state'') asymptotically stable and $y^+$ (``threshold'') unstable. Any initial data $y(0)>y^+$ will traverse the circle before returning to rest.  As $b$ decreases to $b_{SN}=1$, the two roots merge and then for $b<1$, $dy/dt>0$ always and there is a limit cycle.  Henceforth, we will model the excitable cells by (\ref{eq:exc}) with $b>1$.  Oscillators are modeled as the simple phase dynamics,
$$\frac{dx}{dt}=\omega$$ 
where $\omega>0$ is the natural frequency. As with the excitable system, $x\in[0,2\pi)$ lies on the circle.  In general, we will study small chains of excitable systems driven by oscillators at each end and then analyze the locking patterns
\begin{eqnarray}
\label{eq:oeNo}
\begin{aligned}
\frac{dx}{dt} &= \omega+d + c_{oe}\sin(y_1-x) \\
\frac{dy_1}{dt} &= f(y_1) + c_{eo}\sin(x-y_1) + c_{ee}\sin(y_2-y_1) \\
\frac{dy_j}{dt} &= f(y_j)+c_{ee}[\sin(y_{j-1}-y_j)+\sin(y_{j+1}-y_j)] \\
\frac{dy_N}{dt} &= f(y_N) + c_{eo}\sin(z-y_N) + c_{ee}\sin(y_{N-1}-y_N) \\
\frac{dz}{dt} &= \omega-d +c_{oe}\sin(y_N-z)
\end{aligned}
\end{eqnarray}
where $j=2,\ldots, N-1$. Here $x,z$ are oscillators (often referred to as {\bf O cells}) with uncoupled frequencies of $\omega\pm d$ and the variables $y_j$ are excitable (referred to as {\bf E cells}) with $b>1$.  The coupling strength between cells are positive, that is $c_{eo}$, $c_{oe}$, $c_{ee} > 0$. We allow for some heterogeneity in the oscillators via the parameter $d$, also positive. While this may seem as a somewhat restricted parameterization for a model, we note the normal form for a SNIC bifurcation is
\[
\frac{dx}{dt} = 1-\cos(x)+\big(1+\cos(x)\big)p = (1+p)\Big(1-\frac{1-p}{1+p}\cos(x)\Big)
\]
which, after rescaling time, is identical to our model dynamics. For the biophysical simulations, we use the Morris-Lecar model, where each cell obeys
\begin{eqnarray}
\begin{aligned}
V' &= I-4m_\infty(V)(V-120)-8 w (V+84)\\ {} &-2(V+60)+I_{coup}& \\ 
w' &=  0.3\big(w_\infty(V)-w\big)/\tau_w(V) \\
m_\infty(V)&=\frac{1}{2}\Big(1+\tanh\big((V+1.2)/18\big)\Big) \\
w_\infty(V)&=\frac{1}{2}\Big(1+\tanh\big((V-12)/17.4\big)\Big) \\
\tau_w(V)&=\mbox{sech}\big((V-12)/34.8\big)
\end{aligned}
\label{eq:mlmodel}
\end{eqnarray}
with $I=43$ for the oscillators and $I=39$ for the excitable cells.  Coupling currents, $I_{coup}$ have the form $g (\hat{V}-V)$ where $\hat{V}$ is the voltage of the cell to which $V$ is coupled.  The parameter $g$ varies and is provided in the figure captions.  

\section{\label{sec:level1}Results}

Here we outline the results for various numbers of excitable units.  Henceforth, we say that an excitable cell {\em fires} if it traverses the circle passing through $y=\pi$.    

We first explore one OE pair to see the effects of the oscillator on an excitable unit and then look at what happens with chains of excitable cells.

\subsection{\label{sec:level2}OE pair}

We start with the simple system
\begin{eqnarray}
\begin{aligned}
\dot{x} &=  \omega+c_{oe}\sin(y-x) \\
\dot{y} &=  f(y) + c_{eo}\sin(x-y)
\end{aligned}
\label{eq:oe}
\end{eqnarray}
where we set $b=1.1, \omega=1$ and vary the coupling parameters $c_{oe}$ and $c_{eo}$.  This is a system on a two-dimensional torus and as long as $c_{oe}<1$, there are no fixed points.  Since this is a flow on a torus and $\dot{x}>0$, we can make a Poincare section along an arbitrary value $x=C$ which will lead to a one-dimensional map. As the dynamics are in the plane, the map is monotone and invertible, thus, there is a well defined rotation number
\[
\rho = \lim_{t\to\infty} \frac{y(t)}{x(t)}
\] 
which is a continuous function of the parameters. 
When $c_{eo}$ is sufficiently small (e.g., $c_{eo}<(b-1)$ is sufficient), then $y(t)$ will just oscillate around $y^-$, the stable rest state, and the rotation number is 0.  

\begin{figure}
\includegraphics[width=.5\textwidth]{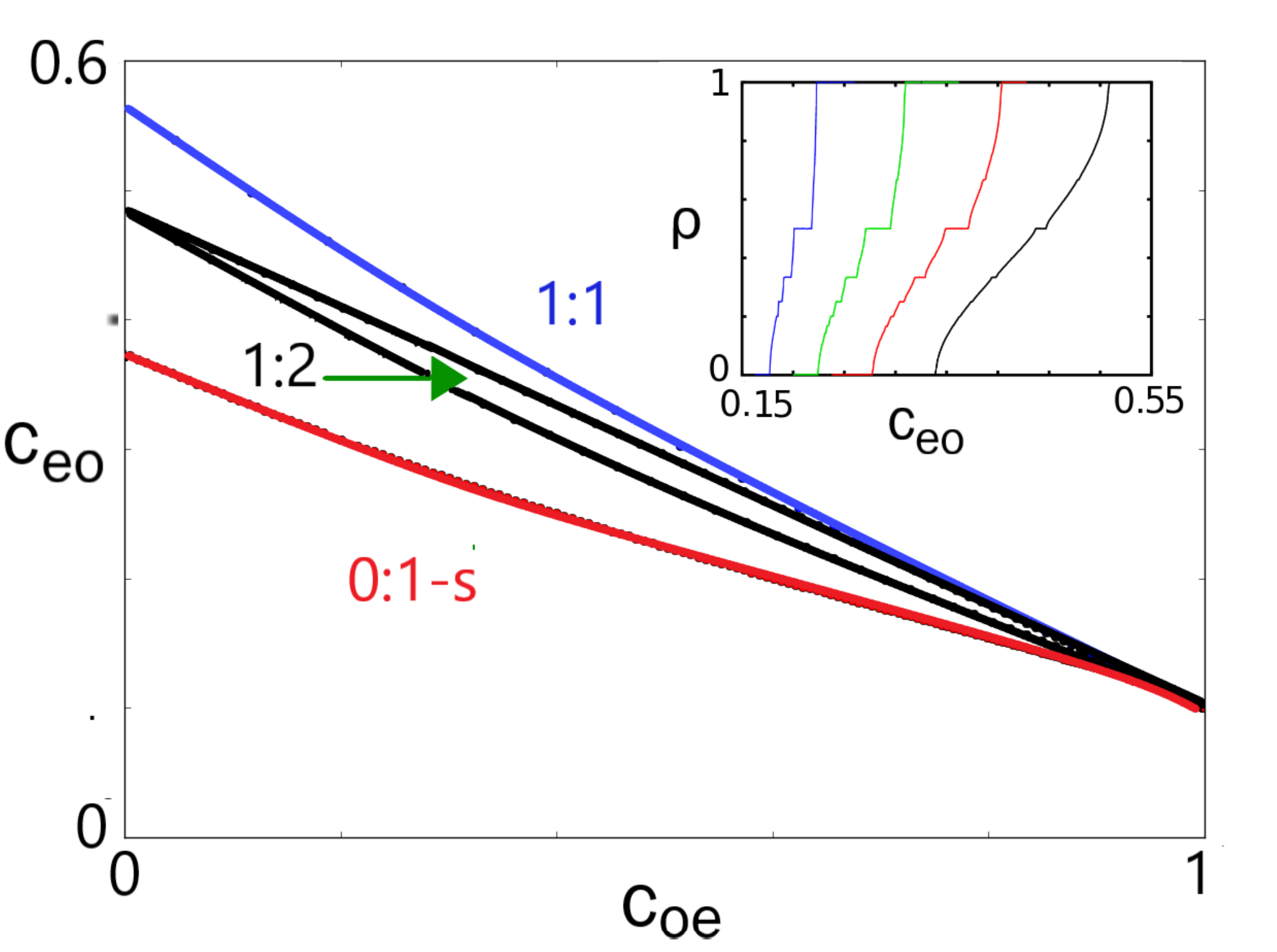}
\caption{Dynamics of Eq. (\ref{eq:oe}) as the connectivity varies. Below the red (bottom) curve, the excitable cell does not fire, while above the blue (top) curve, it fires in a 1:1 manner with the oscillator. In between, rational and irrational firing patterns occur; the 1:2 locking regime is illustrated in between the two black (middle) lines. In fact, all n:m lockings will occur for $n \leq m$, however the regions may be very small. Inset shows the rotation number for different values of $c_{oe}$ as a function of $c_{eo}$. In the picture, $c_{oe}=\{0.1,0.3,0.5,0.7\}$ from left to right.}
\label{fig:oe}
\end{figure}

Figure \ref{fig:oe} shows the behavior of (\ref{eq:oe}) as the coupling strengths vary. If $c_{eo}$ is small enough, then the excitatory cell will never fire, while for $c_{eo}$ large enough, it will always fire in a 1:1 manner with the oscillator.  As $c_{oe}$ goes to 1 (the uncoupled frequency of the oscillator), the oscillator slows its frequency to 0 and in this case, $x$ becomes nearly constant.  With $x$ slowly varying, the equation for $y$ can be treated adiabatically so that $c_{eo}\sin(x-y)$ is a constant lying between $[-c_{eo},c_{eo}]$.  Thus, when $c_{eo}$ exceeds $b-1$, $\dot{y}$ will be positive and the excitable system will fire.  This explains why all the curves in the figure converge at $c_{eo}=b-1$ when $c_{oe}=1.$  The inset in the figure shows the rotation number as a function of $c_{eo}$ at different values of $c_{oe}.$  Higher values of $c_{oe}$ slow down the oscillator so that the critical coupling threshold approaches the minimum value of $b-1.$  At the other extreme, when $c_{oe}=0$, then the rotation number has no open sets of parameters where there are locking regimes other than 0:1 and 1:1.   

We note that if $c_{oe}\ge 1$, then it is possible to find equilibria in equation (\ref{eq:oe}). Multiplying the $\dot{y}$ equation by $c_{oe}$ and the $\dot{x}$ equation by $c_{eo}$ and adding results in
\[
c_{eo}\omega + c_{oe} = c_{oe}b \cos(y).
\]
Thus, fixed points $\bar{y}$ exist as long as $c_{oe}>c_{eo}\omega/(b-1)$. Furthermore, we also must have that $c_{oe}\ge \omega$, since otherwise $\dot{x}>0$. Thus, there is a critical value of $c_{eo}=b-1$ where there is a saddle-node bifurcation with $c_{oe}=\omega$.  In general, the saddle-node bifurcation is $c_{eo}=(b-1)c_{oe}/\omega$ for $c_{oe}>\omega.$  
We close this section by noting that making $b$ larger shifts the curves in Figure \ref{fig:oe} toward higher values of $c_{eo}$ as it takes stronger coupling to induce the excitable cell to fire.  Decreasing the uncoupled frequency of the oscillator from $\omega=1$ is similar to increasing the coupling $c_{oe}$ as both slow the oscillator down giving the excitable system a better chance at firing.

\subsection{\label{sec:level2}OEO chain}

\begin{figure}[h]
\includegraphics[width=.5\textwidth]{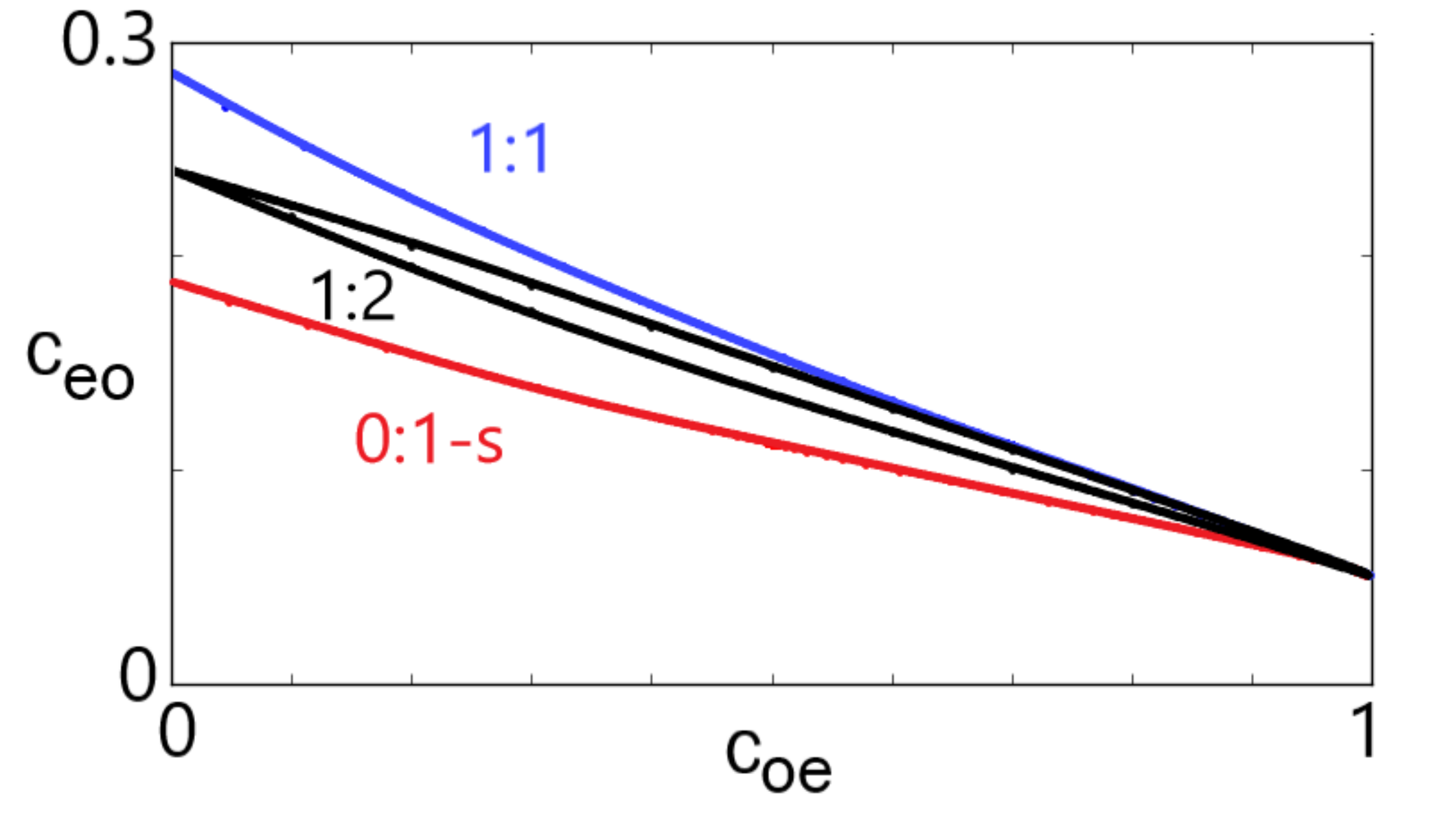}
\caption{Dynamics of Eq. (\ref{eq:oeo}) as the connectivity varies. Again, below the red (bottom) curve, the excitable cell does not fire, and above the blue (top) curve, it fires in a 1:1 manner with the oscillator. In between the black (middle) lines is the 1:2 locking.}
\label{fig:oeobdry}
\end{figure}

The simplest way that two oscillators can interact via an excitable cell is given by
\begin{eqnarray}
\label{eq:oeo}
\begin{aligned}
\dot{x}&= \omega+d + c_{oe}\sin(y_1-x) \\
\dot{y_1} &= f(y_1) + c_{eo}(\sin(x-y_1)+\sin(z-y_1)) \\
 \dot{z}&= \omega-d + c_{oe}\sin(y_1-z). 
\end{aligned}
\end{eqnarray}
We set $b=1.1,\omega=1,d=0$ and varied $c_{oe},c_{eo}$ to get a big picture of the dynamics. Figure \ref{fig:oeobdry} shows boundaries for these phase-locked solutions. For most values of the coupling parameters, the dominant behaviors are 1:1 and 0:1 where the E cells either fire on every cycle or don't fire at all.  Within a narrow sector of parameters, we find the 1:2 phase-locking, where the E cell fires once for every two times the oscillators fire, just like the OE system. We also find 1:3 locking but only in a narrow band of $(c_{oe},c_{eo})$ values; note this existed in the OE system too but it was also too small to show. The boundaries of the OEO system are not very much different from the OE system, although it takes smaller values of $c_{eo}$ for the E cell to fire due to its receiving two oscillatory inputs.

We note that in all choices of $(c_{oe},c_{eo})$ there was always synchrony between the oscillators $x$ and $z$. Thus, asymptotically, with $x(t)=z(t)$, equation (\ref{eq:oeo}) is identical to (\ref{eq:oe}) with $c_{eo}$ doubled. 

\subsubsection{Heterogeneity}

\begin{figure}[h]
\includegraphics[width=3.1in]{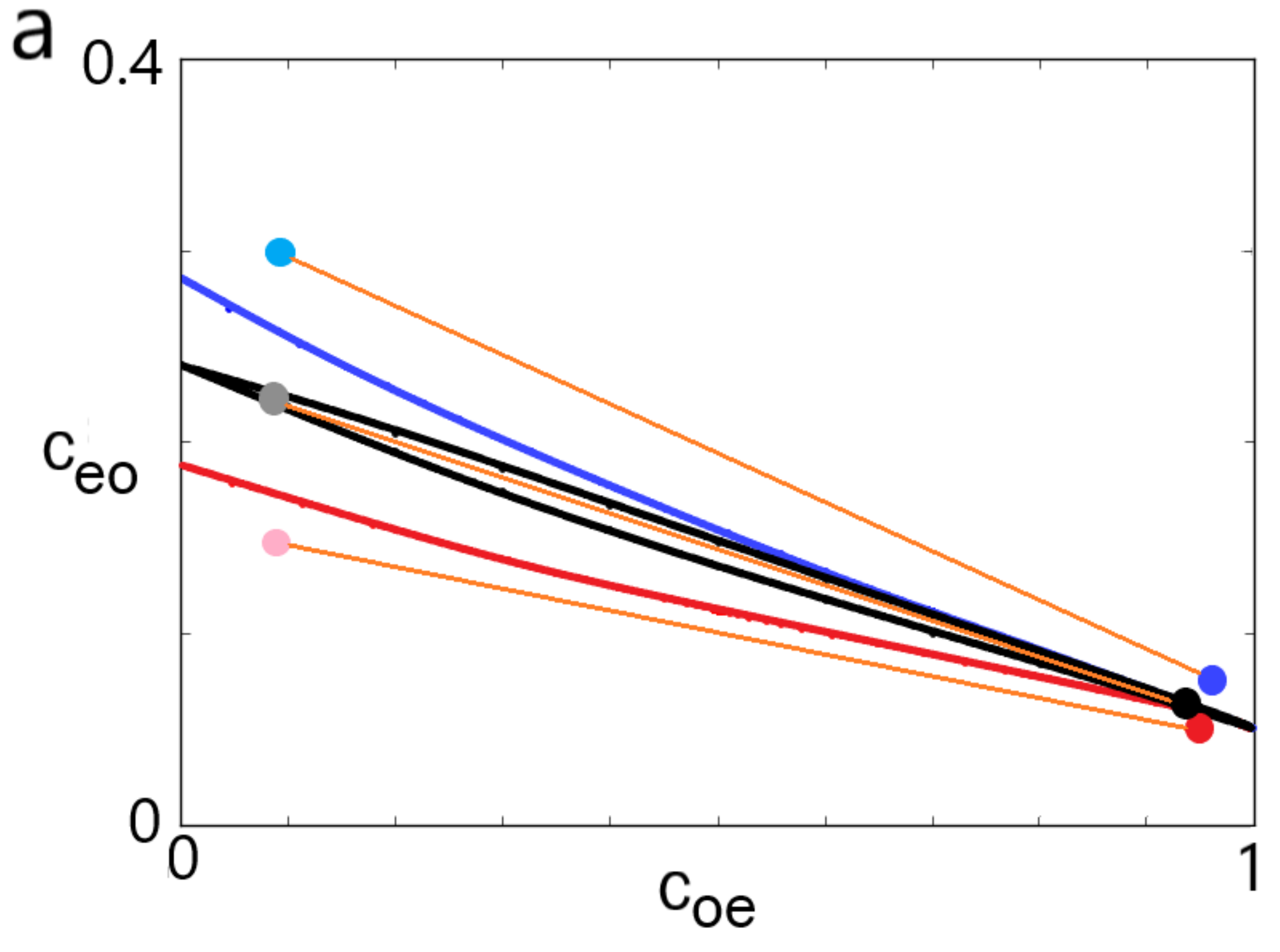}
\includegraphics[width=3.1in]{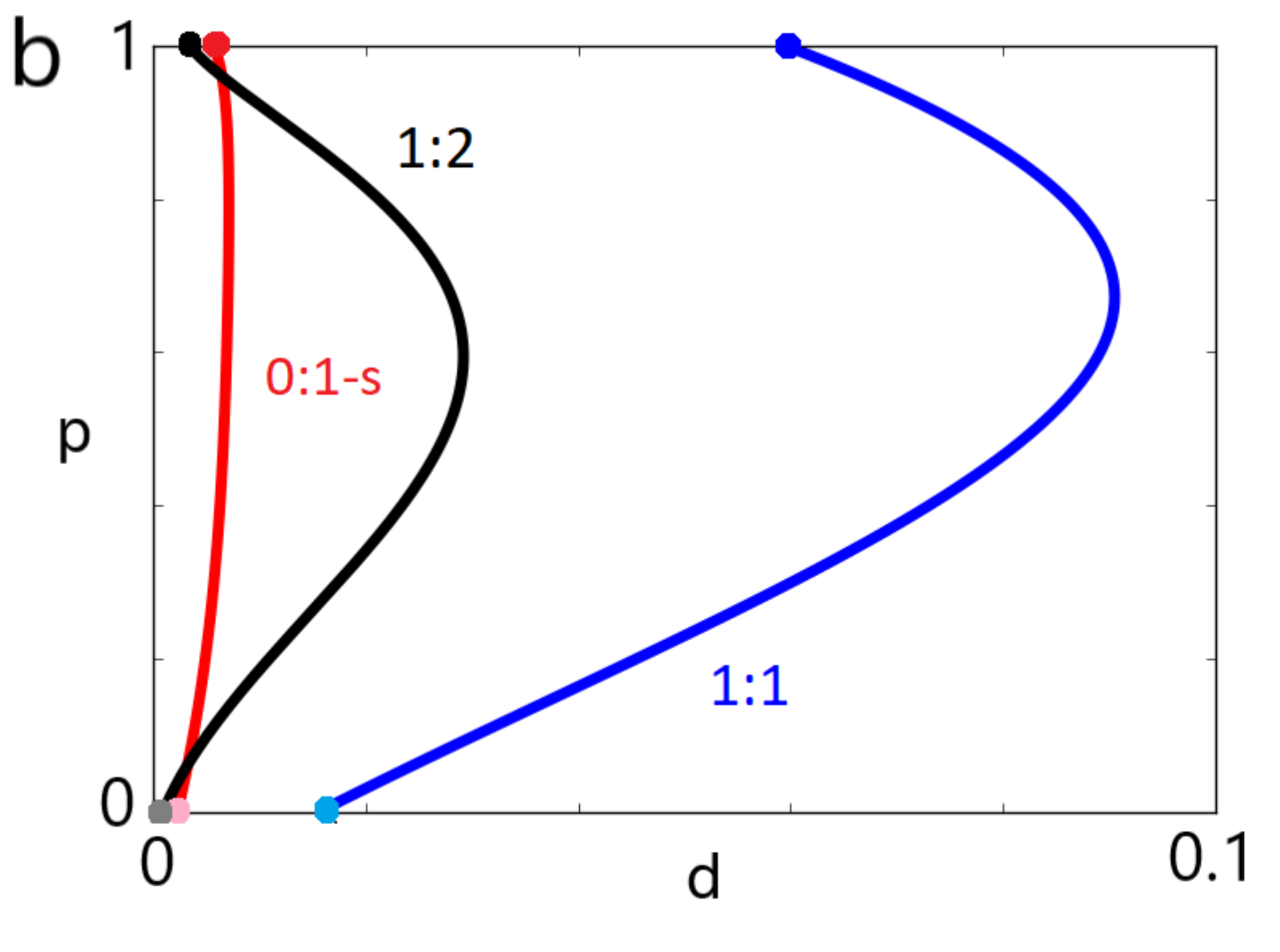}
\caption{The top figure shows the homotopy lines for the OEO chain. These were computed as $(c_{oe},c_{eo})=(k_1,c_1)+p(d_k,d_c)$ for $0\leq p\leq1$. Light colored dots on the left of each line correspond to $p=0$ and dark colored dots on the right of each line are when $p=1$. For the 0:1-s homotopy, the bottom line goes from $(c_{oe},c_{eo})=(0.1,0.15)$ to $(0.95,0.05)$. The line between the black (1:2) goes from $(0.1,0.22)$ to $(0.95,0.06)$. Lastly, the line in the 1:1 region on the top goes from $(0.1,0.3)$ to $(0.95,0.075)$. In the bottom picture, the curves correspond to the boundary when the specified phase-locked solution switches from stable to unstable as $d$ increases. It is clear that 1:1 is more robust than the other solutions and as long as the E cell fires, the locked solution will exist for a wider range of  $d$.}
\label{fig:oeohtp}
\end{figure}

In Figure \ref{fig:oeohtp}, we explore how the change in oscillator frequency affects existence of the phase-locked solutions. Rather than vary $c_{oe}$ or $c_{eo}$, we have chosen to co-vary them along the lines shown in Figure \ref{fig:oeohtp}a as this guarantees that the locking pattern is constant. This also allows us to explore the efficacy of the E cell in coupling the two O cells.  From Figure \ref{fig:oeohtp}b, when the E cell fires, it has a much greater effect on the O cells and thus allows them to lock over a much wider range of heterogeneity. We remark that there is a ``sweet'' spot for coupling strength along this line that maximizes the allowable heterogeneity. Since $c_{oe},c_{eo}$ varies along a straight line, this point is not where the total coupling, $c_{eo}+c_{oe}$, is maximal; that occurs at either $p=0$ or $p=1.$  

\subsection{\label{sec:level2}OEEO chain}

With two E cells, we obtain the equations
\begin{eqnarray}
\label{eq:oeeo}
\begin{aligned}
\dot{x} &= \omega + d+ c_{oe}\sin(y_1-x) \\
\dot{y}_1 &= f(y_1) + c_{ee}\sin(y_2-y_1)+c_{eo}\sin(x-y_1) \\
\dot{y}_2 &= f(y_2) + c_{ee}\sin(y_1-y_2)+c_{eo}\sin(z-y_2) \\
\dot{z} &= \omega -d + c_{oe}\sin(y_2-z).
\end{aligned}
\end{eqnarray} 
There is now one more parameter, $c_{ee}$, which governs the strength of connectivity between the two E cells and thus is important in communicating between the two O cells. As above, we will set $\omega=1$ and restrict all the coupling parameters to lie in $(0,1)$.  

The addition of another E cell makes the dynamics much more complex with multiple stable attractors. If $d=0$, then $x=z$ and $y_1=y_2$ (the synchronous solution) is invariant under the dynamics of equation (\ref{eq:oeeo}) and, in this case, it reduces to the dynamics or equation (\ref{eq:oe}).  However, if this synchrony manifold is unstable, then we can expect to see more complicated behavior.

\begin{figure}[h]
\includegraphics[width=3.5in]{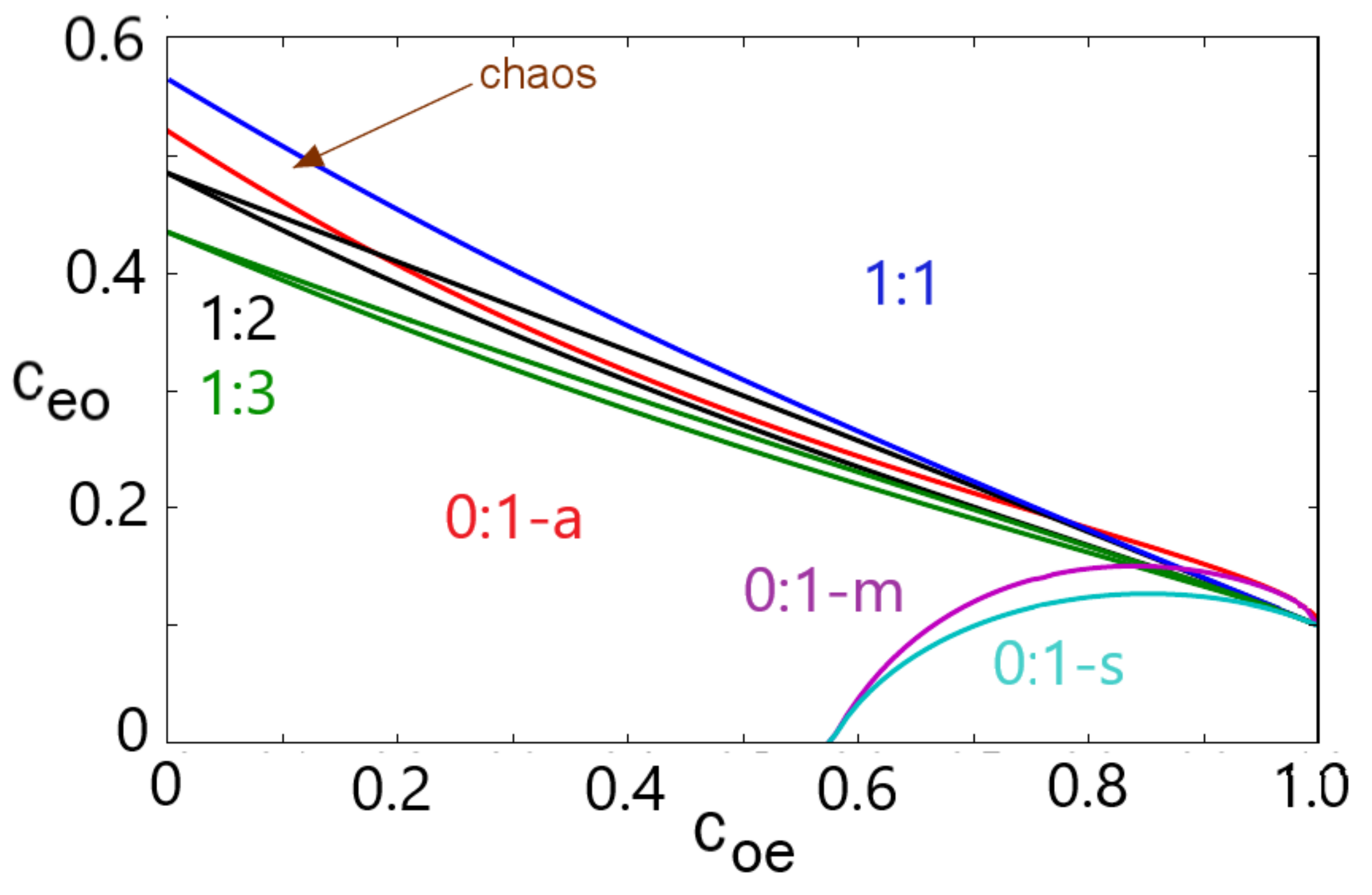}
\caption{Regions of different phase locking for the OEEO model when $c_{ee}=0.5$ as a function of the parameters $c_{eo}$ and $c_{oe}$.  Everything above the blue curve is synchrony with $y_{1,2}$ firing in 1:1 with $x,z$. Between the red curve and the magenta curve, $x,z$ fire in anti-phase and $y_{1,2}$ do not fire. Within the black curves, $x,z$ are synchronous and $y_{1,2}$ fire in a 1:2 manner. Within the green curves, $x,z$ are synchronous and $y_{1,2}$ fire in a 1:3 manner. Between magenta and cyan, $x,z$ have a mixed phase-difference and $y_{1,2}$ don't fire. Finally, below the cyan curve, $x,z$ are synchronous and $y_{1,2}$ don't fire. In other regions such as $c_{oe}=0.11,c_{eo}=0.49$, we have found apparent chaotic behavior. }
\label{fig:oeeo-full}
\end{figure}

Figure \ref{fig:oeeo-full} shows a two-parameter diagram of some of the behavior for $c_{ee}=0.5$ in the $(c_{oe},c_{eo})$ plane. The simplest types of dynamics are $n$:$m$-locking regimes where the E cells fire $n$ cycles for every $m$ cycles of the O cells. For the synchronous dynamics shown in this figure, 0:1,1:1,1:2 and 1:3 locking of the E cells to the O cells all appear to be attractors of this system. The lines in the figure split the $(c_{oe},c_{eo})$ plane into regions of stability. Above the blue line, there is 1:1 synchronous behavior; between the black lines, there is 1:2 synchronous behavior, between the green lines we have 1:3 synchrony, and below the cyan curve, 0:1 synchrony. 

\begin{figure}[h]
\includegraphics[width=3.5in]{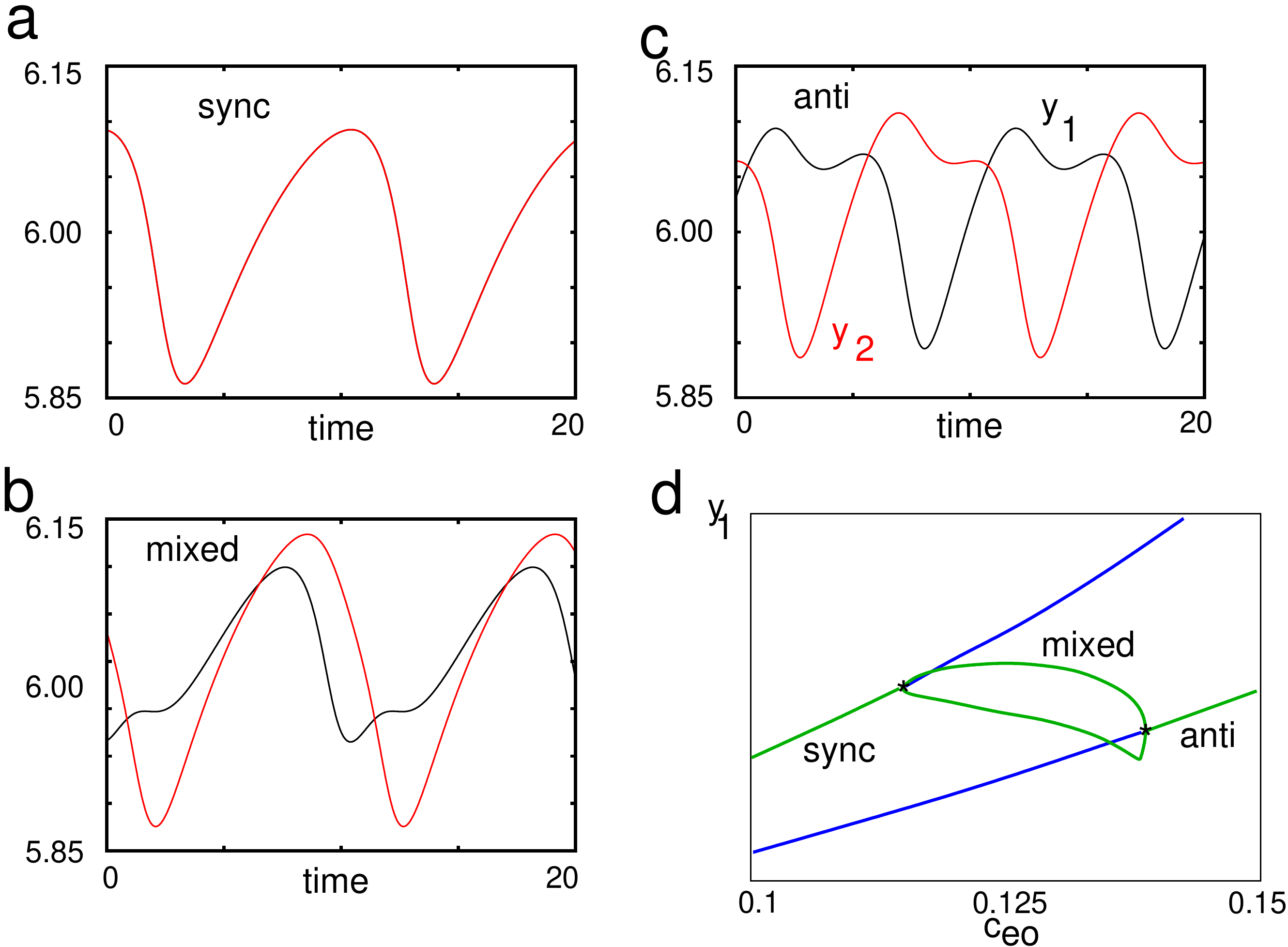}
\caption{The three types of sub-threshold dynamics for $y_{1,2}$ with $c_{oe}=0.78,c_{ee}=0.5$. In (a), $c_{eo}=0.1$, we have 0:1 synchrony. In (b), $c_{eo}=0.13$, we have the ``mixed'' state which is neither synchronous nor anti-phase. In (c) with $c_{eo}=0.15$, we have the anti-phase state. Lastly, (d) shows the pitchfork (symmetry-breaking) bifurcation diagram showing the emergence of the mixed state. The top line in the synchrony branch and as $c_{eo}$ increases, the line changes from stable periodic orbits (green) to unstable periodic orbits (blue). While on the bottom line, the anti-phase branch changes from unstable periodic orbits to stable periodic orbits as $c_{eo}$ increases.}
\label{fig:oeeo-subth}
\end{figure}

\begin{figure}
\includegraphics[width=3.5in]{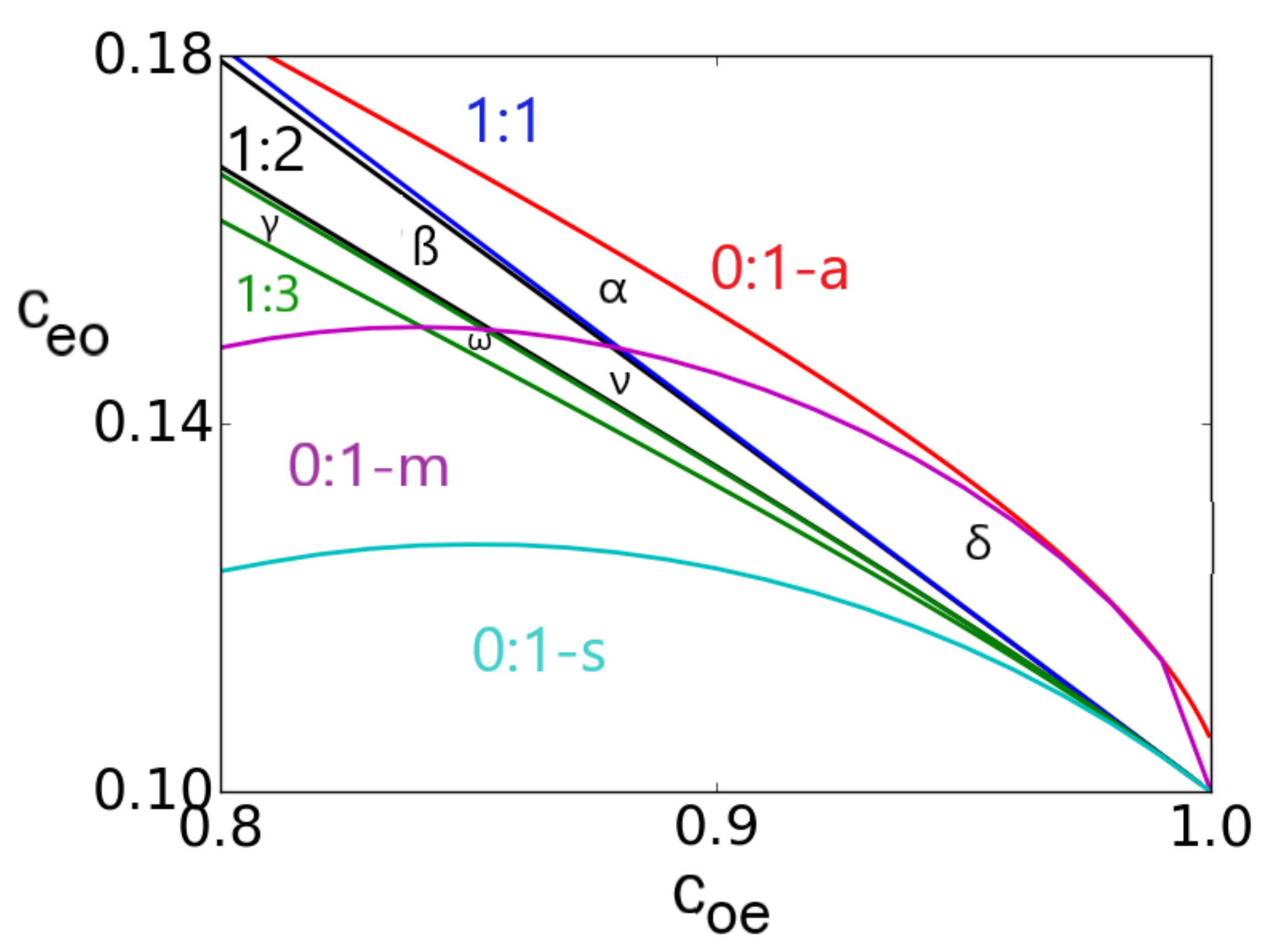}
\caption{Zoomed in picture of Figure \ref{fig:oeeo-full} showing different regions of bistability. In the region labeled $\delta$, there is both 1:1 synchrony and 0:1 mixed, while in region $\alpha$, 1:1 synchrony is bistable with antiphase (0:1-a). Regions $\beta$ and $\gamma$, antiphase behavior coexists with 1:2 and 1:3 locking respectively and regions $\nu$ and $\omega$, mixed sub-threshold behavior coexists with 1:2 and 1:3 locking respectively.}
\label{fig:oeeo-bistable}
\end{figure}

In addition to the synchronous behavior, we also find other stable behavior, where the E cells do not fire. There appear to be three distinct types of this behavior: synchrony (0:1s), anti-phase (0:1a), and ``mixed'' (0:1m).  This is quite different than the OEO system where we were unable to find any stable behavior when $x$ and $z$ weren't synchronized. Anti-phase exists and is stable throughout the region bounded by the red curve and the magenta curve. Synchrony without the E cells firing is stable below the cyan curve. Between the cyan and the magenta curve, we find the so-called ``mixed'' state. Figure \ref{fig:oeeo-subth} shows these three types of behavior. We can best understand the mixed state as follows. Fix $c_{oe}$ at, say $c_{oe}=0.78$ and $c_{eo}$ at a value below the cyan curve where there is stable 0:1 synchrony. Increasing $c_{eo}$ (a vertical line in Fig. \ref{fig:oeeo-full} at $c_{oe}=0.78$) results in a pitchfork or symmetry-breaking bifurcation where a stable branch of non-synchronous asymmetric orbits arises with a phase-difference between synchrony and anti-phase. This is shown in Figure \ref{fig:oeeo-subth}d.  The magenta and cyan curves in Fig \ref{fig:oeeo-full} depict these pitchfork bifurcations. (A zoomed in version near $c_{oe}=1$ is shown in Figure \ref{fig:oeeo-bistable}.)  

\begin{figure}[h!]
\includegraphics[width=3.5in]{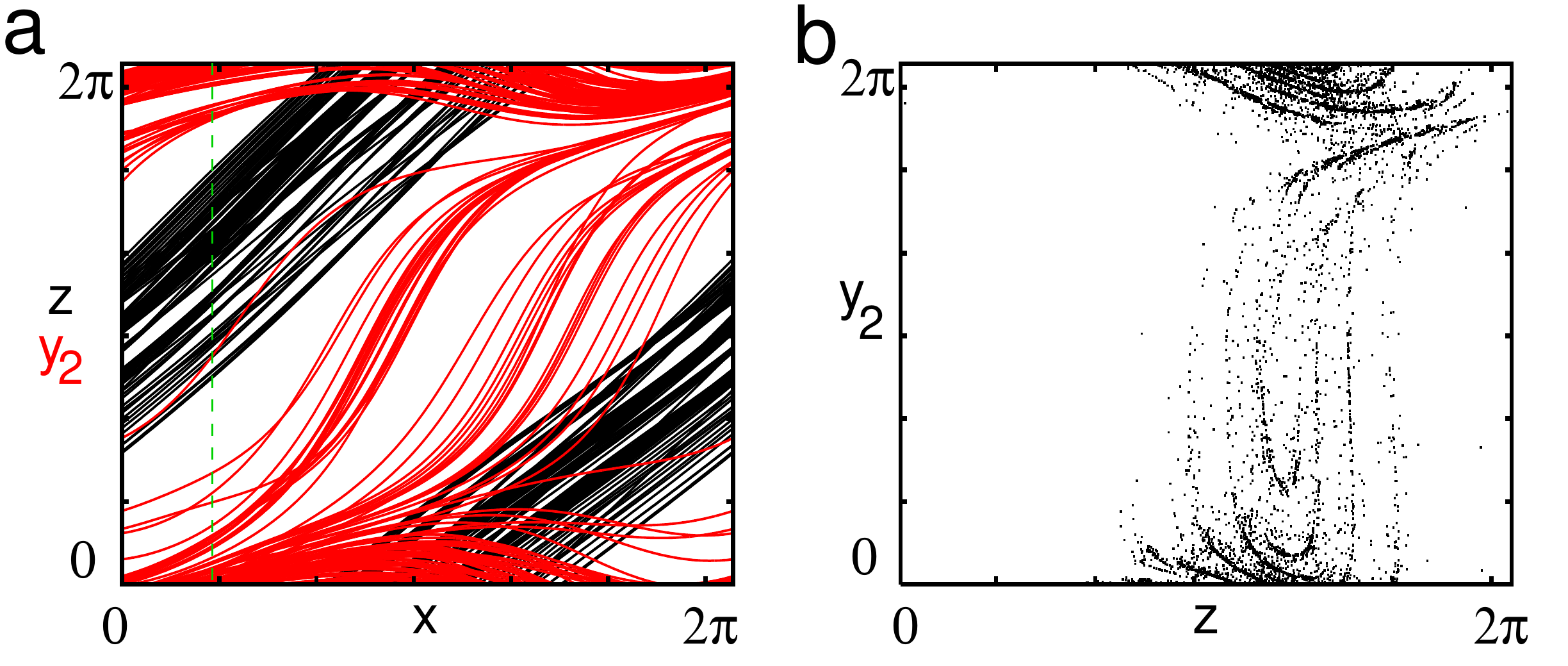}
\caption{Putative chaos in Eq.~(\ref{eq:oeeo}) when $c_{oe}=0.11,c_{eo}=0.49,c_{ee}=0.5$. Picture (a) shows phase space trajectories of $z$ and $y_2$ vs $x$. Picture (b) is a Poincare section through $x=1$, the thin dotted green line in panel (a), showing $y_2$ vs $z$.}
\label{fig:oeeo-chaos}
\end{figure}

Other regions not accounted for include the region below the 1:1 synchrony line and above anti-phase line and the top 1:2 line. This large space does not have any apparent phase-locked pattern and appears to be chaos. For example, when $c_{ee}=0.5,c_{oe}=0.11,c_{eo}=0.49$, Figure \ref{fig:oeeo-chaos} shows $z$ and $y_2$ vs $x$ in panel (a) and a Poincare section through $x=1$ in panel (b). We have crudely estimated the Liapunov exponent to be about 0.04. If we restrict Eq. (\ref{eq:oeeo}) to the synchrony manifold, then we find quasi-periodic behavior for $y_1(t)$. Our general observation is that stable $1:m$ locking for Eq. (\ref{eq:oe}) leads to stable synchrony for Eq. (\ref{eq:oeeo}).

\subsubsection{Bistability}

The OEEO chain is the smallest chain we have found that exhibits regions of bistability, that is, the long term dynamics depend on the initial data. Figure \ref{fig:oeeo-bistable} shows a zoomed in version of Figure \ref{fig:oeeo-full}. Regions labeled by Greek letters indicate regions of bistability. Since 1:1 synchrony is stable above the blue line and 0:1 antiphase  (0:1-a) is stable below the red curve and above the magenta curve, we see that in the region labeled $\alpha$, there is bistability between these two states.  Similarly, in region $\beta$ (respectively, $\gamma$), both 0:1-a and 1:2 (resp. 1:3) are stable. In region $\nu$ (respectively $\omega$), the 0:1 mixed (0:1-m) and 1:2 (resp. 1:3) are stable. 

\begin{figure}[h!]
\includegraphics[width=1.6in]{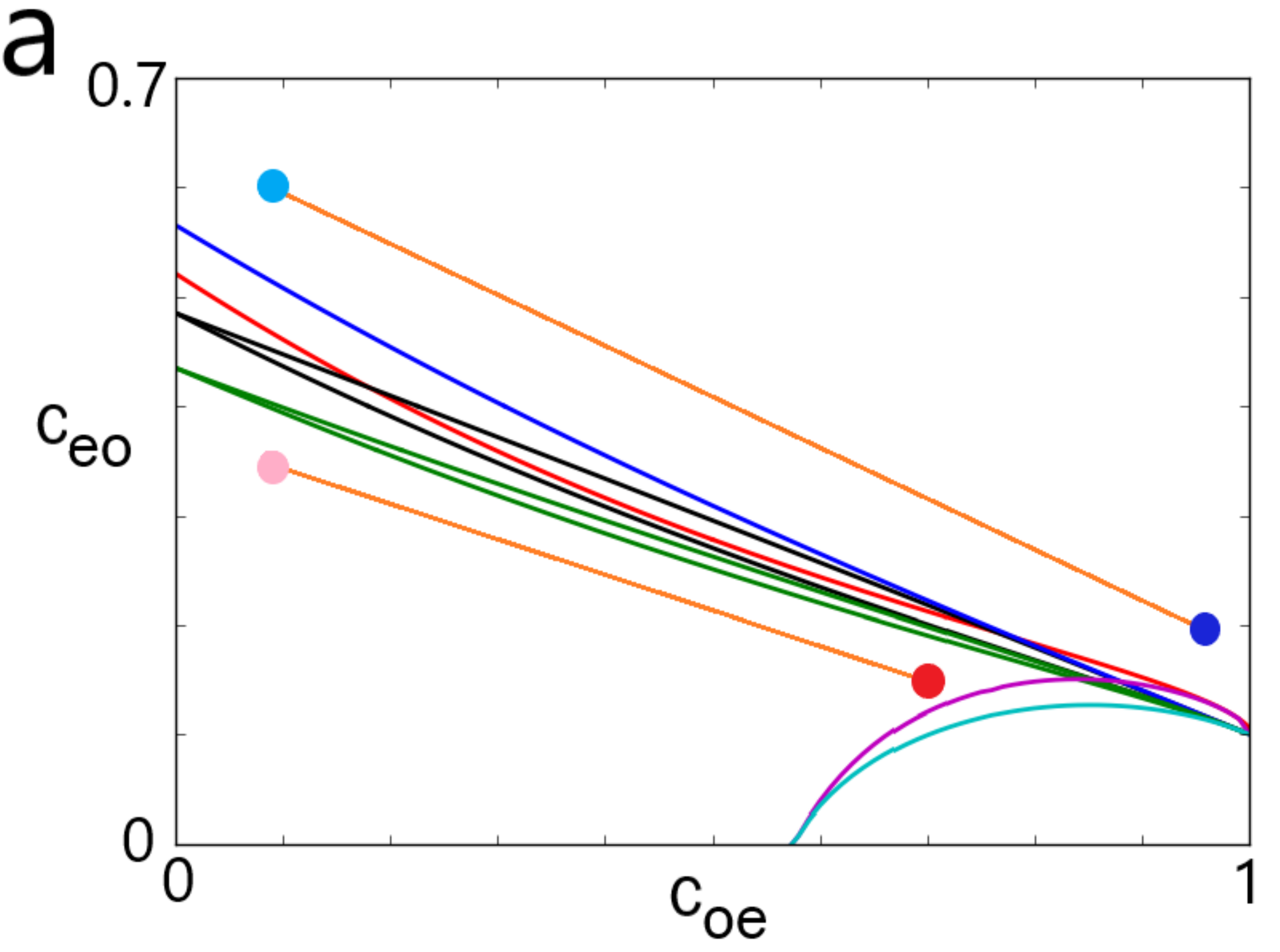}
\includegraphics[width=1.6in]{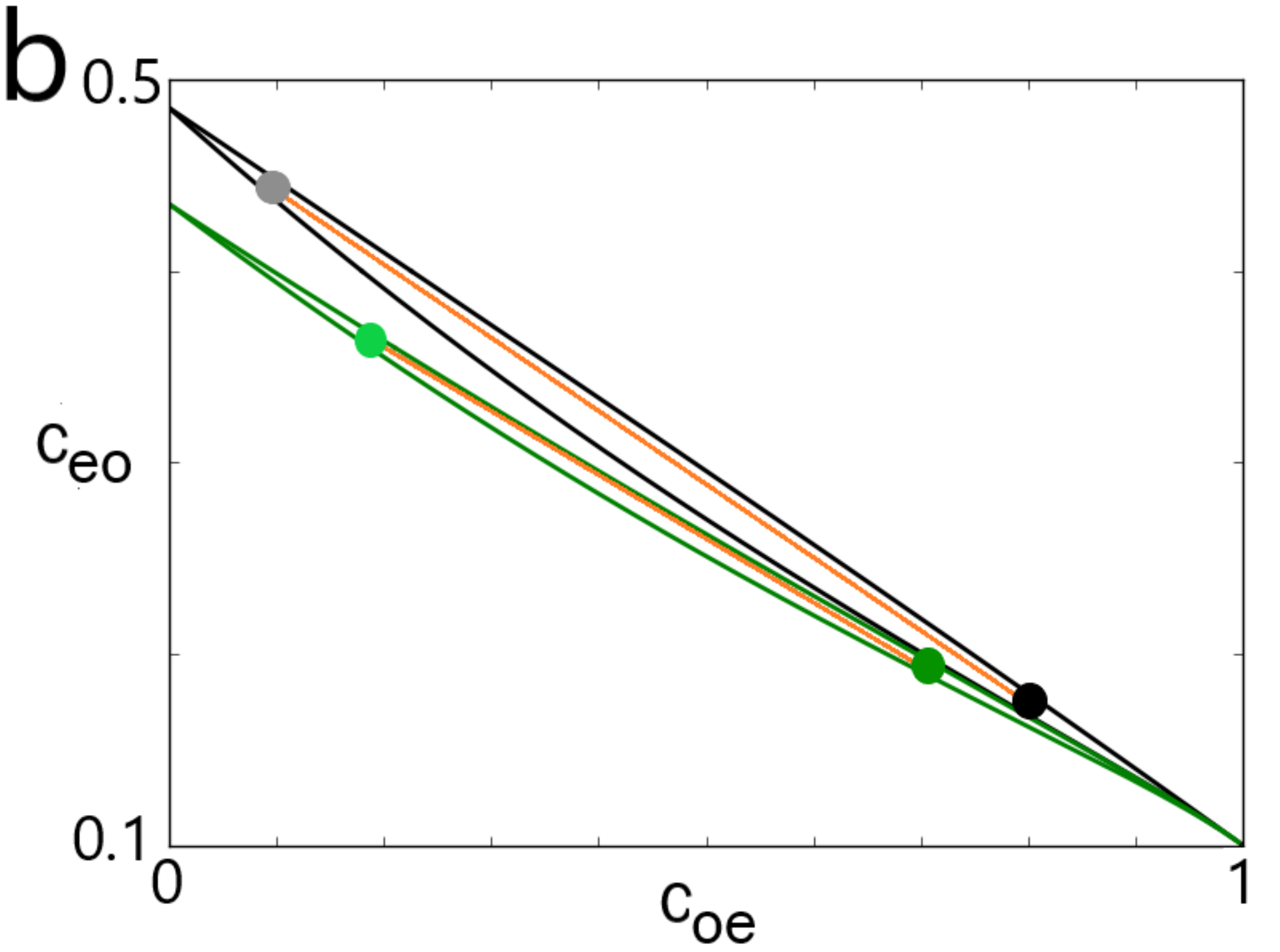}
\includegraphics[width=3.3in]{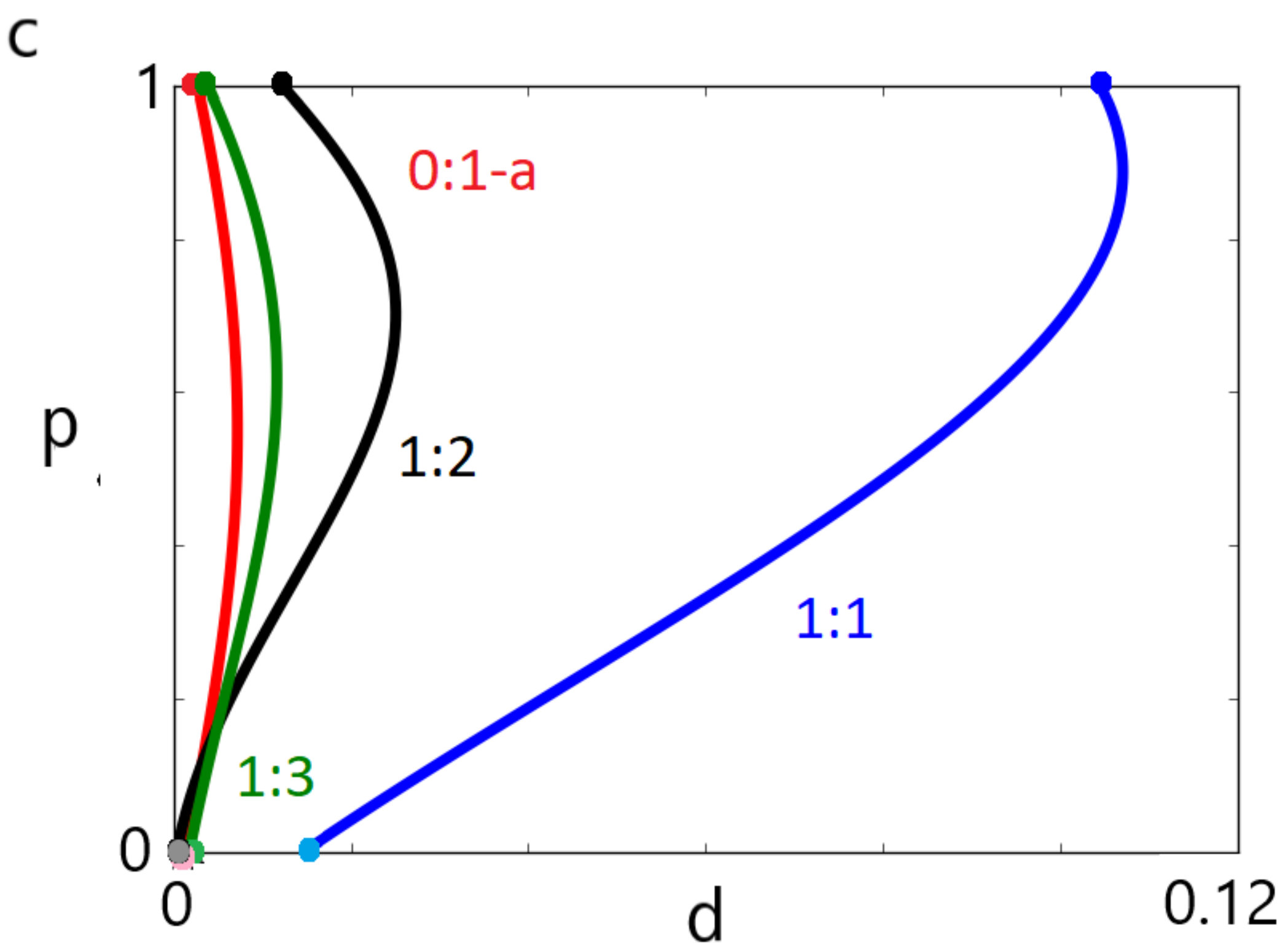}
\caption{In these pictures, we show where we take the homotopy for the non-bistable regions. For the 1:1 locking, our parameters go from $(c_{oe},c_{eo}) = (0.1,0.6)$ to $(0.95,0.2)$ and for 0:1-a locking, our parameters go from $(0.1,0.35)$ to $(0.7,0.15)$. The second picture shows our homotopies for the 1:2 and 1:3 regions. For these, we have our parameters going from $(0.1,0.44)$ to $(0.8,0.174)$ and $(0.2,0.36)$ to $(0.71,0.19)$, respectively. Here we show how large the heterogeneity can be in each of the non-bistable regions. The dots on the bottom of this figure correspond to the dots labelling the left endpoint of the lines in (a) and (b) and the dots on the top of this figure correspond to the dots labelling the right endpoint of those lined. As our homotopy parameter $p$ increases, we increase $c_{oe}$ and decrease $c_{eo}$ according to the orange diagonal lines in (a) and (b). We can see that 1:1 coupling is still the most robust and 0:1-a is the least robust despite it having a large stability region.}
\label{fig:oeeo-htp}
\end{figure}

\begin{figure}[h!]
\includegraphics[width=3.3in]{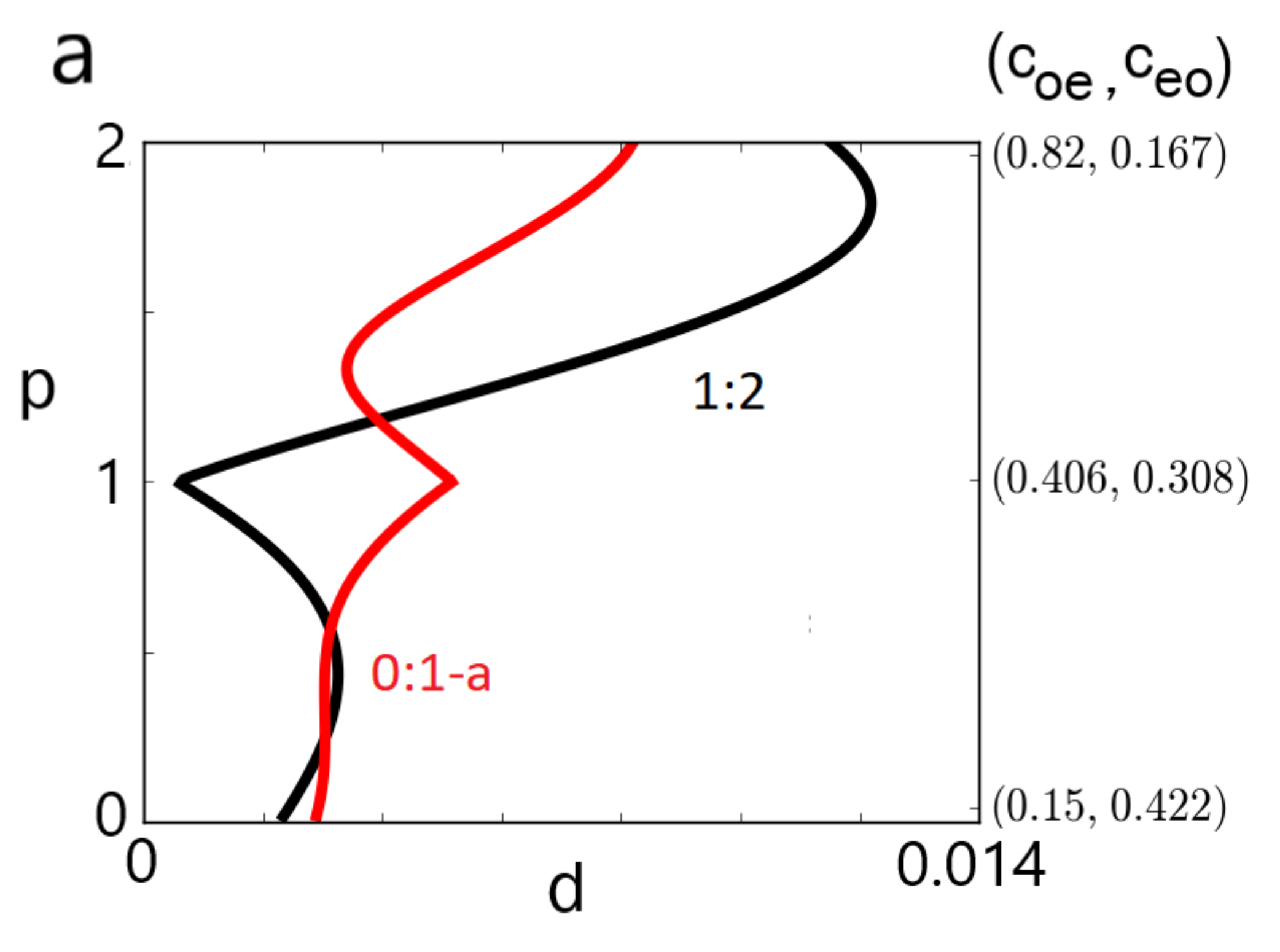}
\includegraphics[width=3.3in]{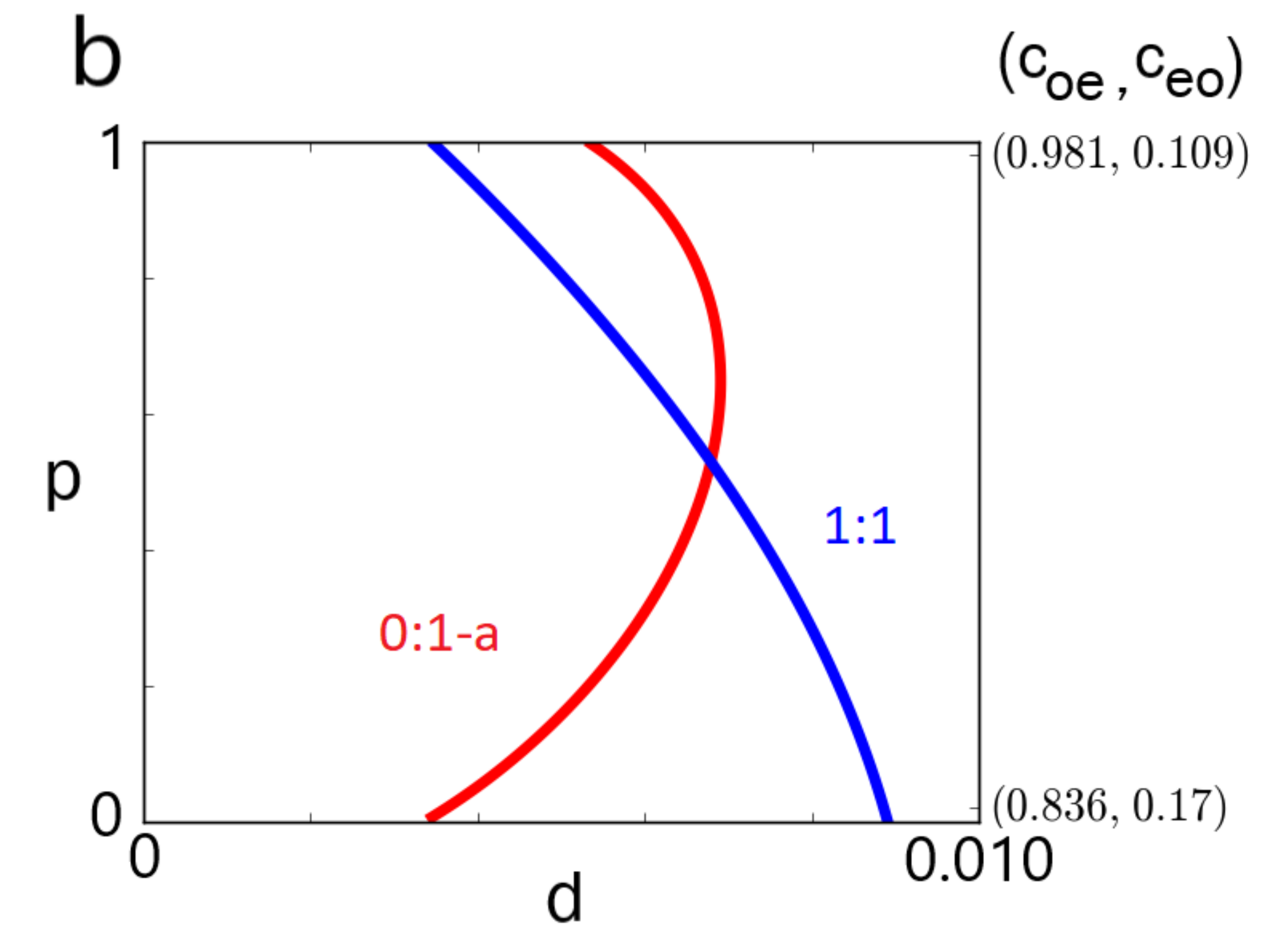}
\caption{Effects of heterogeneity in bistable regions. The right $y$-axis shows the homotopy starting and ending points for $(c_{oe},c_{eo})$. In (a), the first homotopy begins at $(0.15, 0.422)$ in between the black lines and below the red line in Figure \ref{fig:oeeo-full} where 0:1-a and 1:2 are bistable, and goes to $(0.406, 0.308)$, where it is very narrow but still bistable. The second homotopy starts at this point and continues into the $\beta$ region (see Figure \ref{fig:oeeo-bistable}) where it ends at $(0.82,0.167)$. The black 1:2 line (resp. red 0:1-a line) shows the maximum $d$ can be before losing stability of the 1:2 (resp. 0:1-a) locking. In (b), we perform a homotopy in the $\alpha$ region (see Figure \ref{fig:oeeo-bistable}), where 0:1-a and 1:1 are bistable. When $d\neq0$, the 0:1-m region deflates and thus, the $\alpha$ and $\delta$ region become one. Our homotopy went from $(0.836,0.17)$ to $(0.981,0.109)$. The blue 1:1 line (resp. red 0:1-a line) represents the maximum heterogeneity before 1:1 (resp. 0:1-a) stability is lost. It is interesting to note that the curves in both panels cross nontrivially. Parameters $b$ and $c_{ee}$ remain constant at $b=1.1$ and $c_{ee}=0.5$.}
\label{fig:oeeoH-bi}
\end{figure}

\subsubsection{Heterogeneity}

We can apply similar numerical analyses to the OEEO chain as with the OEO chain and compare the existence of locked solutions as the oscillator frequencies vary from $\omega$ in (\ref{eq:oeeo}). However, unlike the OEO chain, we can also investigate bistable regions as well; this will give us a more direct comparison since in these regions, the parameters can remain the same. First, we look at the 4 long term dynamics separately. Figure \ref{fig:oeeo-htp}a and \ref{fig:oeeo-htp}b shows the four homotopies we take for each of the four main regions: 1:1, 1:2, 1:3, and 0:1-a. Using these, Figure \ref{fig:oeeo-htp}c shows the frequency behavior of the four homotopies. It is clear that the 1:1 synchrony region is the most robust based on changes in the frequency. What is interesting is although the 1:2 and 1:3 regions are much smaller than the 0:1-a region, there is more tolerance  with 1:2 and 1:3 as the heterogeneity parameter increases. This could suggest that the excitable cells firing plays a key role in maintaining locking as $d$ increases. This is similar to what we saw in the OEO system with Figure \ref{fig:oeohtp}. As with the OEO system, there is a ``sweet'' spot where the system is most tolerant of frequency differences.

We can also look at the regions of bistability. There are two regions we looked into: when 0:1-a and 1:2 were both stable and when 0:1-a and 1:1 were both stable. Besides the $\beta$ region in Figure \ref{fig:oeeo-bistable}, 1:2 and 0:1-a are also bistable for $c_{oe}$ smaller (see Figure \ref{fig:oeeo-full} or \ref{fig:oeeo-htp}a). Furthermore, these regions are the same region, one can see in Figure \ref{fig:oeeo-full} that the 0:1-a and 1:2 bistable region does not break as $c_{oe}$ increases. So we did two homotopies for this long skinny region: one line for the top half and one line for the bottom half and we made sure the two lines connected. Figure \ref{fig:oeeoH-bi}a shows what happens. As the homotopy moves towards the very narrow region in the middle (near $c_{oe}$=0.4), the 1:2 becomes less tolerant with $d$ and the 0:1-a region becomes more tolerant. Then, as the homotopy enters the $\beta$ region, 1:2 locking allows a bigger range of  $d$ before becoming unstable. The other region of bistability we looked into was the $\alpha$ region in Figure \ref{fig:oeeo-bistable}. We can see initially the 1:1 synchrony and more tolerable but as we move down into what was the $\delta$ region, it is the 0:1-a dynamics that are more robust as $d$ increases (see Figure \ref{fig:oeeoH-bi}b).

\subsection{\label{sec:level2}Weak coupling via E cells}

\begin{figure}[h]
\includegraphics[width=2.9in]{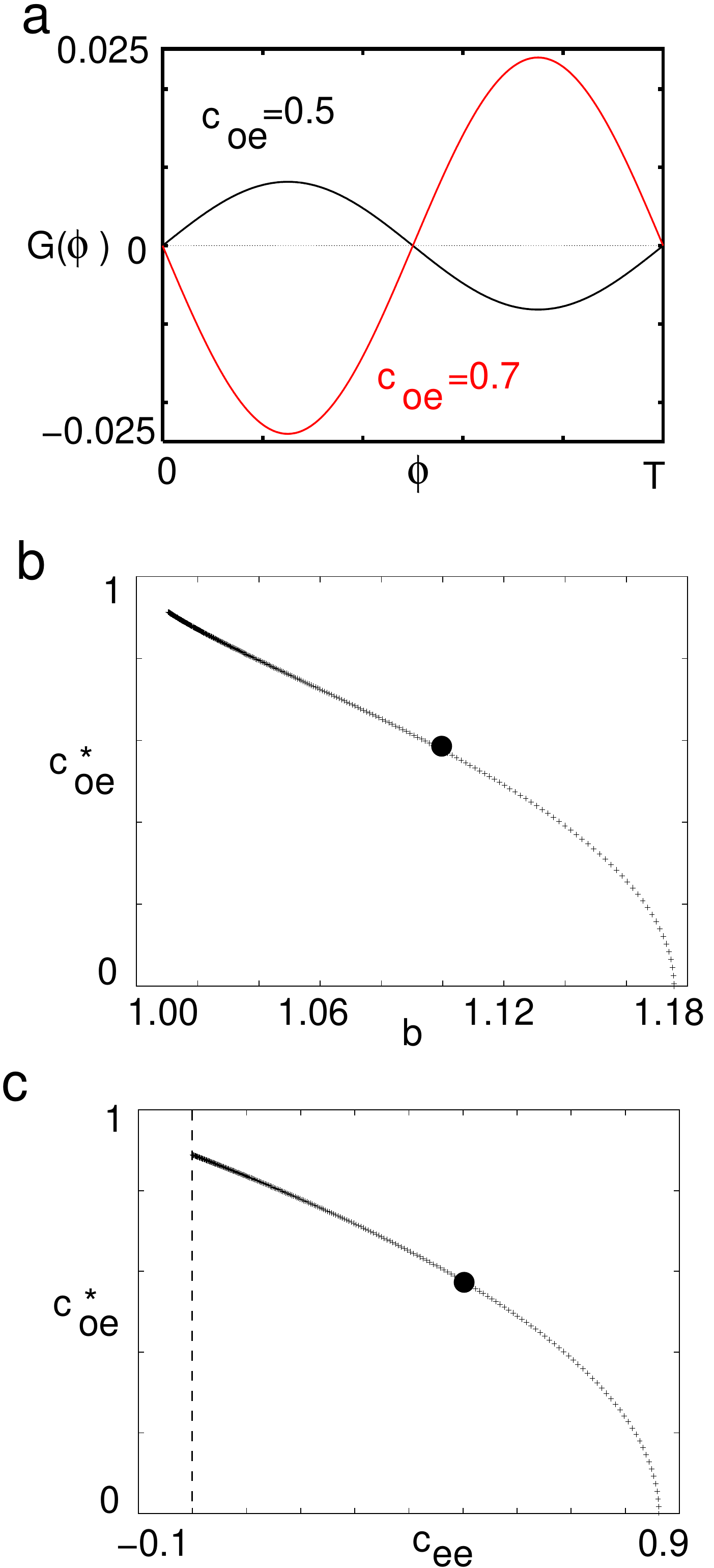}
\caption{The first picture shows $G(\phi)$ as a function of $\phi$. The different graphs describe the flipped behavior for $c_{oe}=0.5$ and $c_{oe}=0.7$. In picture B and C, we see how the critical $c_{oe}$ value varies with $c_{ee}$ and $b$. In this section, $b$ has been at 1.1 and $c_{ee}$ has been at 0.5. For $c_{oe}>c_{oe}^{*}$ (resp $c_{oe}<c_{oe}^{*}$), synchrony is stable (unstable) and anti-phase is unstable (stable). The black dots in the bottom two graphs signify the parameters we have used: $b=1.1$ and $c_{ee}=0.5$. }
\label{fig:oeeoweak}
\end{figure}

When $c_{eo}$ is small enough, the E cells do not fire and there can be several types of dynamics including synchrony, anti-phase, and a non-synchronous locked state (see Figure \ref{fig:oeeo-full} near $c_{oe}=0.57$).  Based on this figure, it appears that in the limit as $c_{eo}\to0$, there is an abrupt transition from anti-phase to synchrony as $c_{oe}$ crosses a critical value.  We will now address this point using weak coupling analysis. This method extends to arbitrary length chains but for simplicity, we just perform the analysis for the OEEO chain. Let $c_{eo}=\epsilon$ where $0<\epsilon\ll 1$ is a small parameter and so we have
\begin{eqnarray*}
x' &=& 1 + c_{oe}\sin(y_1-x)    \\
z' &=& 1 + c_{oe}\sin(y_2-z)  \\
y_1' &=& f(y_1) + \epsilon \sin(x-y_1) + c_{ee}\sin(y_2-y_1) \\
y_2' &=& f(y_2) + \epsilon \sin(z-y_2) + c_{ee}\sin(y_1-y_2),
\end{eqnarray*}    
and as a reminder, $f(q)=1-b\cos(q)$. When $\epsilon=0$, we can set $y_{1}=y_{2}$ to be constant and we have $y_{1,2}\equiv k = -\arccos(1/b)$.  If $0<c_{oe}<1$, there is a $T-$periodic solution, $U(t)$, to $U'=1-c_{oe}\sin(U)$ with $U(t+T)=U(t)+2\pi.$ Note that 
\[
T=\int_0^{2\pi} \frac{dx}{1-c_{oe}\sin(x)}.
\]
Thus for $\epsilon$ small, we expect that $x(t)\approx k + U(t+\theta_x)$ where $\theta_x$ is an arbitrary phase shift. To formalize this argument, we use the method of multiple scales by  letting $s=t$ be the fast time, and $\tau=\epsilon \hspace{1pt} t$ be the slow time and expand $x,z,y_{1,2}$ as a power series in $\epsilon$, e.g., 
$$x(t)=x_0(s,\tau) + \epsilon x_1(s,\tau) + \ldots,$$ 
$$y_1(t)=y_1^0(s,\tau) + \epsilon y_1^1(s,\tau)+\ldots,$$ 
where we assert each term is $T-$periodic in $s$.  With this ansatz, in the first order expansion we see that $x_0(s,\tau)=k + U(s+\theta_x(\tau))$, $z_0(s,\tau)=k+U(s+\theta_z(\tau))$, and $y_{1,2}=k$ where $\theta_{x,z}(\tau)$ are unknown. In the second order expansion,

\begin{eqnarray}
\label{eq:wx}
\partial_s x_1 + U'(s+\theta_x) \partial_\tau \theta_x =c_{oe}\cos(U(s+\theta_x))[y_1^1-x_1] \nonumber \\
\partial_s z_1 + U'(s+\theta_z) \partial_\tau \theta_z =c_{oe}\cos(U(s+\theta_z))[y_2^1-z_1] \nonumber \\
\label{eq:wy}
\partial_s y_1^1 = b\sin(k) y_1^1 + c_{ee}(y_2^1-y_1^1) +\sin(U(s+\theta_x)) \nonumber \\
\partial_s y_2^1 = b\sin(k) y_2^1 + c_{ee}(y_1^1-y_2^1) +\sin(U(s+\theta_z)) \nonumber,
\end{eqnarray}
where $\displaystyle \partial_\gamma = \frac{\partial}{\partial \gamma}$. The last two equations can be written as
\[
\frac{\partial}{\partial s} \left(\begin{array}{c} y_1^1 \\~\\ y_2^1
  \end{array}
  \right) = A\left(\begin{array}{c} y_1^1 \\~\\ y_2^1
  \end{array}
  \right)
 +  \left(\begin{array}{c} \sin\big(U(s+\theta_x)\big) \\~\\ \sin\big(U(s+\theta_z)\big)
  \end{array}
  \right)
\]
where

\[
A = \left(\begin{array}{cc} b\sin(k)-c_{ee} & c_{ee} \\~\\
c_{ee} & b\sin(k)-c_{ee}
\end{array}\right).
\]
This matrix $A$ has strictly negative eigenvalues and thus there is a unique periodic solution to this linear system. Let $W=(w_1,w_2)^T$ be the periodic solution to:
\[
\frac{\partial W}{\partial s} = A W + (\sin(U(s)),0)^T.
\]
Then $y_1^1(s,\tau)=w_1(s+\theta_x)+w_2(s+\theta_z)$ and $y_2^1(s,\tau)=w_1(s+\theta_z)+w_2(s+\theta_x)$.  Now that we have solved for $y_j^1$, we turn to $x_1,z_1$.    
 Consider the linear operator on the space of differentiable $T-$periodic functions:
\[
M(s) x := \partial_s x + c_{oe}\cos(U(s))x.
\]
Since $U'(s)=1-c_{oe}\sin(U(s))$, we see that $x=U'(s)$ is in the
nullspace of $M$ and thus $M$ has a one-dimensional nullspace. With the standard $L^{2}$ inner product, $(f,g) = \int_{0}^{T} f(s)g(s) \hspace{3pt} ds$, the operator $M(s)$ has an adjoint, $M^{*}(s) x = -\partial_s x + c_{oe}\cos(U(s)) x$ and a nullspace, $1/U'(s).$  With this notation, the equation for $x_1(s,\tau)$ can be written as:
\begin{multline*}
M(s+\theta_x) x_1 + U'(s+\theta_x)\partial_{\tau}\theta_{x} \\=c_{oe}\cos(U(s+\theta_x)) [w_1(s+\theta_x)+w_2(s+\theta_z)],
\end{multline*}    
\begin{multline*}
M(s+\theta_z) z_1 + U'(s+\theta_z)\partial_{\tau}\theta_{z} \\=c_{oe}\cos(U(s+\theta_z)) [w_1(s+\theta_z)+w_2(s+\theta_x)].
\end{multline*}  
Taking the inner product of both sides of the $x_{1}$ equation with $1/U'(s+\theta_x)$, we obtain the dynamics of $\theta_x$:
$$
T \partial_\tau \theta_x = \int_{0}^{T} \frac{c_{oe}\cos(U(s+\theta_x))}{U'(s+\theta_x)} [w_1(s+\theta_x)+w_2(s+\theta_z)]\ ds.
$$
A simple change of variables gives $\partial_\tau \theta_x = H(\theta_z-\theta_x)$ where 
\begin{equation}
\label{eq:hweak}
H(\phi) = \frac{c_{oe}}{T} \int_0^T \frac{\cos(U(s))}{U'(s)} [w_1(s)+w_2(s+\phi)]\ ds.
\end{equation}
Similarly, $\partial_\tau \theta_z = H(\theta_x-\theta_z)$. Finally, we let $\phi=\theta_z-\theta_x$ and use $\partial_\tau \theta_z$ and $\partial_\tau \theta_x$ to obtain the weak coupling equation:
\begin{equation}
\label{eq:phi}
\partial_s \phi = H(-\phi)-H(\phi)=:G(\phi),
\end{equation} 
where $-G(\phi)/2$ is the odd part of $H(\phi).$  In Figure \ref{fig:oeeoweak}a, we plot $G(\phi)$ for $c_{oe}=0.5$ and $c_{oe}=0.7$.  As can be seen from the figure, when $c_{oe}=0.5$, synchrony ($\phi=0$) is unstable and anti-phase ($\phi=T/2$) is stable and the reverse is true for $c_{oe}=0.7$. 

From Figure \ref{fig:oeeo-full}, it appears that the synchrony and anti-phase boundary meet at exactly the same point on the $c_{oe}$ axis. What this means is that synchrony and anti-phase change stability at the same point for this choice of parameters. Stability of synchrony (resp. anti-phase) is lost when $G'(0)$ (resp. $G'(T/2)$) changes from  negative to positive. Denote $m(s)=\cos(U(s))/U'(s)$ so from the definition of $G(\phi)$:
\[
G'(0)=-\frac{2c_{oe}}{T}\int_0^T m(s)w_2'(s)\ ds
\]
\[
G'(T/2)=-\frac{2c_{oe}}{T}\int_0^T m(s)w_2'(s-T/2)\ ds.
\]
Changing $s$ to $s-T/2$ in the second integral and using the observation $m(s+T/2)=-m(s)$ due to the symmetry of the function $\sin(x)$, this shows that $G'(T/2)=-G'(0)$, so that synchrony and anti-phase swap their stability at the same value of $c_{oe}$ independent of any other parameters. Hence in Figure \ref{fig:oeeoweak}b and \ref{fig:oeeoweak}c, we also show how the critical $c_{oe}$-value varies as $b$ or $c_{ee}$ change. If we change the coupling function between the O and the E cells to some more general odd periodic function, say $\sin(x)+a \sin(2x)$, then the symmetry of $m(s)$ is gone and the branches for synchrony and anti-phase will not meet at a point as $c_{eo}\to 0.$ 

\subsection{\label{sec:level2}OEEEO and beyond}

\begin{figure}[h!]
\includegraphics[width=3.6in]{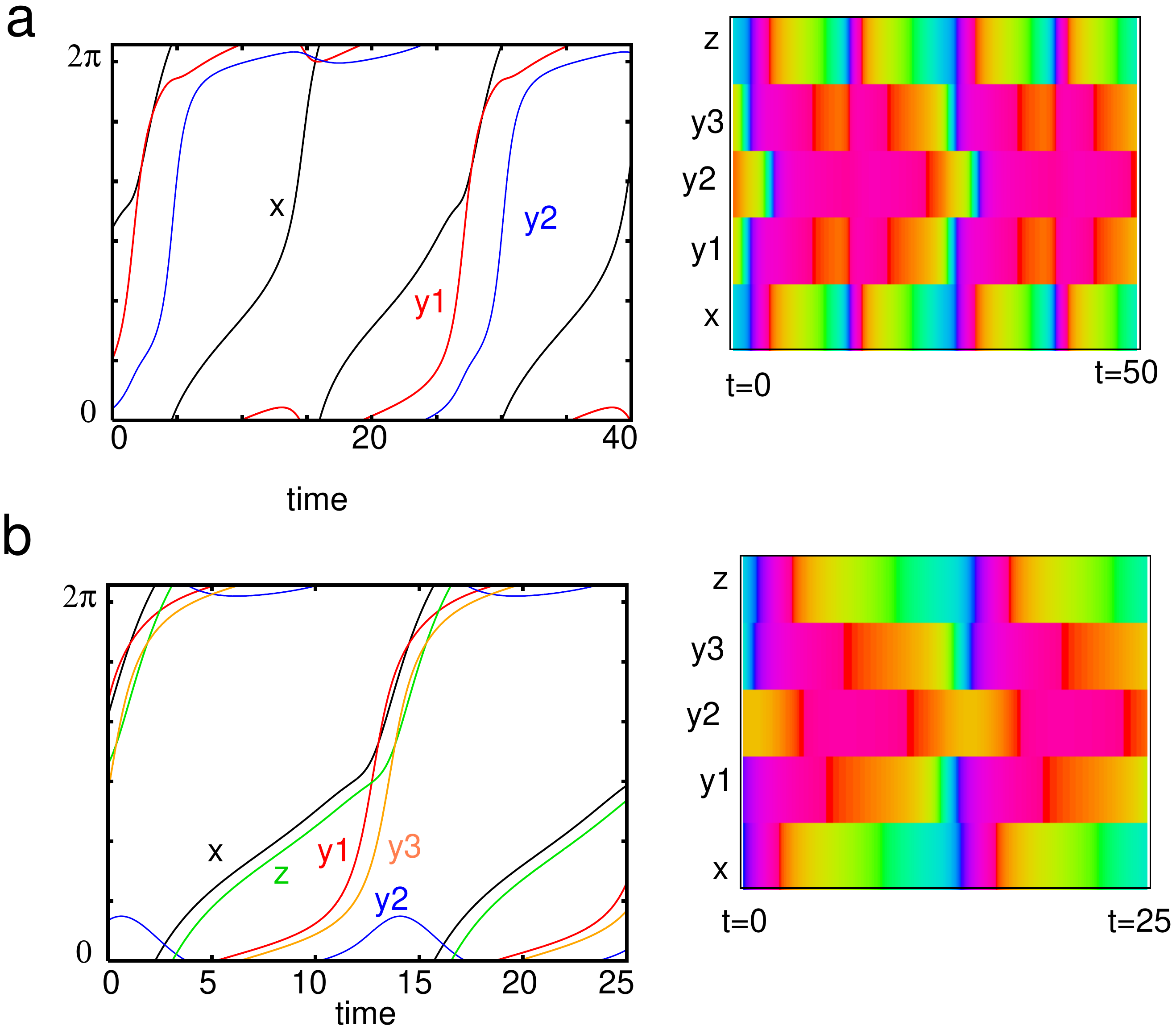}
\caption{Example of bistability in an OEEEO chain. The parameters are $c_{oe}=0.75, c_{eo}=0.25$, and $c_{ee}=0.18$. In (a), $x=z$, and $y_1=y_2$ and $y_m$ (the middle excitable cell) also fires all in synchronous 1:2. In (b), the outer $y_{1,2}$ fire with the O cells in 1:1 but the middle E cell $y_m$ does not fire and the oscillators do not synchronize.}
\label{fig:oeeeobi}
\end{figure}

For chains with 3 or more E cells between the O cells, it is possible to have only the E cells that are coupled to the O cells fire, while the E cells in the middle of the chain fail to fire. For example, Figure \ref{fig:oeeeobi}a shows an OEEEO system where the middle E cell fires and gives rise to 1:2 locking between the oscillators and the excitable cells, while in panel (b), with the same parameters the middle cell $y_m$ does not fire and the O cells fire in 1:1 with the outermost E cells. 

So far, we have seen that for small chains where all the E cells fire, the only stable solution is synchrony and it is robust to small changes in the relative frequencies of the O cells. Let us now consider a long chain of E cells terminated by two O cells acting as pacemakers.  Consider the isolated E chain with no oscillators. If we suppose that the coupling between the E cells is strong enough, then one expects that initiating the first E cell in the chain to fire will result in a traveling wave that propagates down the chain.  If at some time after the first E cell is excited, we initiate a wave at the other end, we expect the two waves to collide somewhere in the middle of the chain and could annihilate completely. This means that the last E cell and first E cell are ``unaware'' the other E cell fired. If we put the two oscillators on the ends, it seems to imply that they will not synchronize; rather they can maintain any phase-difference.  For example, consider
\begin{eqnarray*}
x' &=& \omega_x + c_{oe}\sin(y_1-x) \\
y_1' &=& f(y_1) + c_{eo}\sin(x-y_1) + c_{ee}\sin(y_2-y_1) \\
y_{j}' &=& f(y_j) + c_{ee}\sin(y_{j-1}-y_j) + c_{ee}\sin(y_{j+1}-y_j) \\
y_{100}' &=& f(y_{100}) + c_{eo}\sin(z-y_{100}) + c_{ee}\sin(y_{99}-y_{100}) \\
z' &=& \omega_z + c_{oe}\sin(y_{100}-z)
\end{eqnarray*}

\begin{figure}
\includegraphics[width=3.3in]{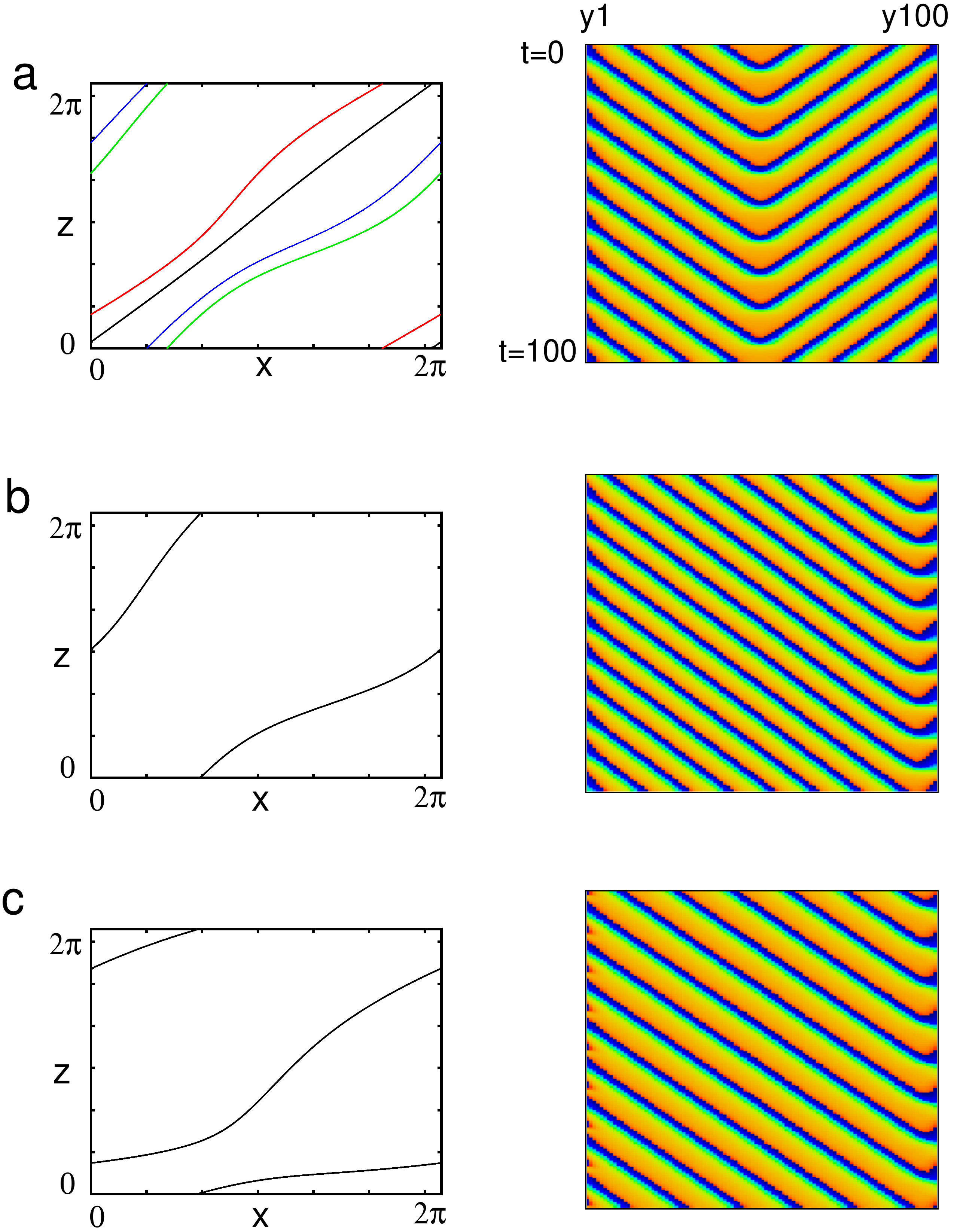}
\caption{Chain of 2 oscillators with 100 excitable cells in between. (a) Phase plane showing the long-time behavior of the two oscillators for 4 different initial data; right is a space-time plot. (b) Same as (a) but $\omega_x=1.1,\omega_z=0.9$; a fixed 1:1 locking always occurs. (c) Same as (a), but $\omega_x=1.5,\omega_z=0.5$ and a 1:2 locking occurs.}
\label{fig:long}
\end{figure} 

\noindent for $j=2,\ldots, 99$. For this section, we took $c_{oe}=0.7, c_{eo}=2, c_{ee}=3$ and $b=1.1$. Figure \ref{fig:long}a shows a simulation when $\omega_x=\omega_z=1$ for four different initial conditions. The left panel shows us that the two oscillators lock but the phase-difference between $x$ and $z$ varies each time. The right panel shows why this happens: waves initiated at the end points collide in the middle and, thus, cannot alter the timing of their opposite oscillators.  However, if we make one oscillator faster than the other, then the point of intersection of the waves moves toward the slower oscillator as the fast oscillator dictates the frequency and becomes a pacemaker. Figure \ref{fig:long}b shows this with $\omega_x=1.1$ and $\omega_z=0.9$. Once there is a single phase-locked solution, there appears to be a unique attractor. Increasing the frequency difference further (Figure \ref{fig:long}c) leads to 1:2 locking where $x$ goes 2 cycles and the rest of the medium goes 1 cycle. Differences in the frequencies of the oscillators allow for the timing information to propagate down the chain and lock the oscillators.   

\subsection{\label{sec:level2}Biophysical Models}

In this paper, we have used a one-dimensional model for excitability that is equivalent to the normal form for a general system near a saddle-node infinite cycle bifurcation (SNIC). A simple and well-known neural model that has a SNIC is the Morris-Lecar (ML) model given in the Methods section. Thus, we turn our attention to this model and look at the OEEO system. Similar to equation (\ref{eq:oeeo}), we hold $c_{ee}=0.1$ and vary $(c_{oe},c_{eo})$ to compare this ML model to the dynamics of equation (\ref{eq:oeeo}).

\begin{figure}[h!]
\includegraphics[width=3.5in]{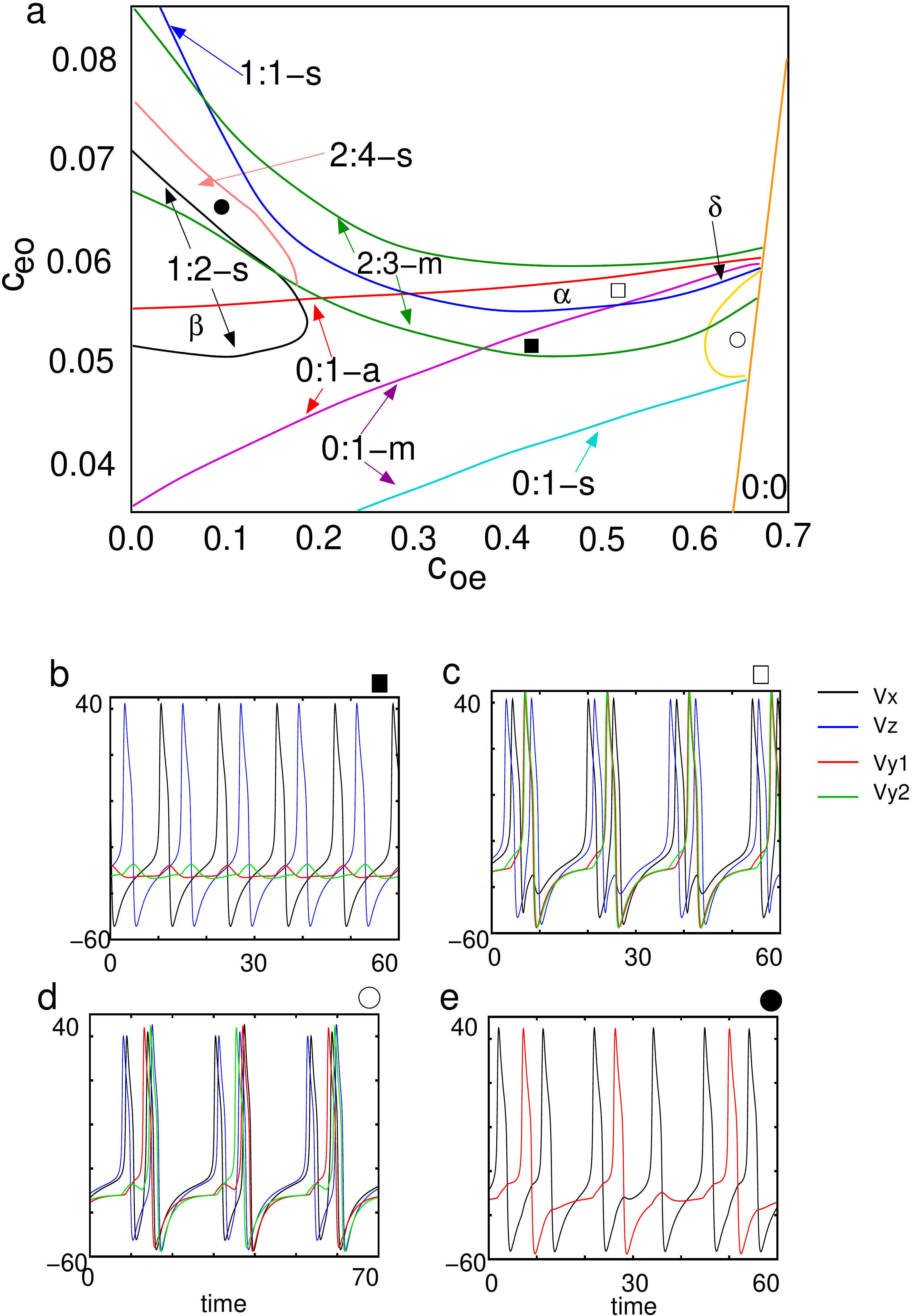}
\caption{(a) Regions of stability for the ML model as constructed from a bifurcation analysis. Stable locking regions are delineated by colored lines with arrows pointing to the boundaries. The symbols $s,m,a$ correspond to the oscillators being synchronous, mixed, or anti-phase. There are many regions of multistability; three are marked $\alpha,\beta,\delta$ corresponding to the regions in Fig. \ref{fig:oeeo-bistable}. Markers correspond to parameters for the time series shown in b-e. Values of $(c_{oe},c_{eo})$ correspond to the markers in (a). (b) 0:1-m,  $(0.4,0.05)$; (c) 2:3-m, $(0.5,0.057)$; (d) 2:4-m, $(0.63,0.05)$; (e)  2:4-s, $(0.1, 0.065)$. In (e), $V_{x}=V_{z}$ and $V_{y_{1}}=V_{y_{2}}$ so they overlap.}
\label{fig:mlfull}
\end{figure}

Figure \ref{fig:mlfull}a shows the regions of stability as we vary $(c_{oe},c_{eo})$ for the ML model. This figure was created by following bifurcation points using AUTO in XPP (see \cite{xpp}), and then combining the two-parameter data and tracing the curves using splines. We show a number of different regions, but this is by no means exhaustive. We compare this figure to Figures \ref{fig:oeeo-full} and \ref{fig:oeeo-bistable}.  We first note that as in the simple phase model, the largest regions correspond to 1:1 synchronous locking and 0:1 locking in either synchrony, mixed, or anti-phase for increasing values of $c_{eo}$. A notable difference from the phase model is the large region of 2:3-m in the ML system. Interestingly, the O cells do not synchronize, but operate in the mixed phase mode. Another difference is that the region of 1:2-s is somewhat limited in the ML model when compared to the phase model.  We have labeled three different regions, $\alpha,\beta,\delta$ in which there is bistability, similar to the phase model in Figure \ref{fig:oeeo-bistable}.  For example in region $\beta$ there is bistability between 1:2-s and 0:1-a.  Both regions $\alpha,\delta$ have bistability between 1:1-s and 0:1. However, due to the existence of the large 2:3-m region, there is actually tristability with the 2:3-m state. The 2:4-s state occurs via a period doubling bifurcation of the 1:2-s state (the upper curve in the 1:2-s region). The 0:1-m state also loses stability via a period-doubling bifurcation as $c_{oe}$ increases (shown by the gold C-shaped region on the right side of the 0:1-m region).  Like the phase model, all of the locked oscillatory regions terminate as a stable fixed point emerges when $c_{oe}$ is large enough. This region is labeled 0:0. Figures \ref{fig:mlfull}b-e show representative voltage traces in some of the different regions indicated by the markers in panel (a).

In summary, we have seen that the simple phase model for interacting oscillatory and excitable cells is a good predictor of the qualitative dynamics of biophysical networks of coupled oscillator and excitable cells. In particular, the latter undergo many of the same bifurcations and transitions between states as well as having similar regions of bistability.

\section{\label{sec:level1}Discussion}

Throughout this project, we have studied some simple networks in which a pair of oscillators is indirectly coupled to active nonlinear elements, namely excitable systems. We show several distinct qualitative behaviors that include in-phase, anti-phase, and mixed-phase synchronization both when the excitable cells fired and did not fire. We also found some regimes of seemingly chaotic dynamics in between phase-locked regions. With the smallest of chains, bistability between the phase-locked regions was impossible; however, if the chain increases in length, many bistable regions can appear. We found that when the excitable cells were silent, that is they operate in the subthreshold regime, the interactions between the two oscillators could be analyzed through weak coupling analysis and was amenable to averaging methods and phase reduction. We found that when the excitable cells are active, the ability to phase-lock is more robust to changes in the oscillator frequencies than when the excitable cells are silent. We also saw that for very long chains of excitable cells, a small difference in the intrinsic frequency of the oscillators is more conducive to rapid and stable phase-locking than if the oscillators are identical.  While this may seem counterintuitive, one can regard it as a case when the faster oscillator becomes the ``leader'' and thus the excitable cells and the slower oscillator are effectively forced. The slow oscillator gets overpowered in the rhythm.  We also showed that the simple phase models that formed the bulk of the paper behaved quite similarly to systems of coupled Morris-Lecar models, a simple biophysical model for a neuron.  

There are many unanswered questions that remain in this paper concerning other types of indirect coupling. In the present paper, we looked at one-dimensional chains. However, a more biologically realistic scenario would involve a small number of oscillators embedded in a two-dimensional network of excitable cells. Indeed, this is a geometry more akin to the examples that motivated this work in the Introduction.

\textbf{Acknowledgement.} This work is supported by NSF grant DMS-1712922.

\end{document}